\documentclass[10pt]{article}
\usepackage[utf8]{inputenc}
\usepackage{latexsym,amssymb,amsmath,url,authblk}

\def\N{\mathbb{N}}

\def\Z{\mathbb{Z}}
\def\R{\mathbb{R}}

\def\proof{\par\noindent{\em Proof. }}
\def\eproof{\hfill{$\Box$}\bigskip}

\def\ds{\dots}
\def\sus{\subset}

\def\cc{\colon}

\newtheorem{thm}{Theorem}[section]
\newtheorem{prop}[thm]{Proposition}
\newtheorem{cor}[thm]{Corollary}

\newtheorem{defi}[thm]{Definition}

\title{Polynomials as terms and the Boolean Independence Theorem}
\author{M. Klazar}
\date{\today}

\begin{document}

\maketitle
\begin{abstract}
We develop a~theory of formal multivariate polynomials over commutative rings by treating them as ring 
terms. Our main result is that two ring terms are s-equivalent (when expanded they 
yield the same standard polynomial) iff they are f-equivalent (one can be 
transformed in the other by a~series of elementary transformations). We consider 
in a~similar way Boolean terms (formulas) and prove a~theorem that two events $a$ and 
$b$ in a~probability space, which are built by two Boolean terms from respective tuples $A$ 
and $B$ of elementary events, are independent if the events in $A$ are 
independent of the events in $B$. This theorem rigorizes arguments in the Probabilistic Method in Combinatorics.  
\end{abstract}

\tableofcontents

\section{Introduction}\label{sec_intro}

We want to rigorize uses of independence of events in the Probabilistic Method in Combinatorics, 
and later we get to it, but it turns out that this task can be joined with and 
motivated by a~better understanding of multivariate polynomials over commutative 
rings, to which we devote the initial larger part of our article.

Let us recall the status of univariate real polynomials  in courses of Mathematical Analysis. 
A~polynomial is a~map $p\cc\R\to\R$ that is either constantly $0$ or is computed by the expression 
$$
{\textstyle
p=p(x)=\sum_{j=0}^n a_jx^j
}
$$
in which $n\ge0$ is an integer, $a_j\in\R$, $a_n\ne0$ and $x$ ranges in $\R$. One denotes the set of all such maps by 
$\R[x]$. We call this definition of $\R[x]$ the MA-{\em definition}. It is widely used but it is also problematic because addition and multiplication in this $\R[x]$ are far from straightforward.

We give a~better analytical definition of $\R[x]$ and call it the EF-{\em definition}. 
For $c\in\R$ let $f_c\cc\R\to\{c\}$ be constant functions
and $\mathrm{id}\cc\R\to\R$, $\mathrm{id}(x)=x$, be the identity 
function. We say that a~map $f\cc\R\to\R$ is a~{\em polynomial} if it can be obtained 
by finitely many applications of two rules.
\begin{enumerate}
    \item $f=f_c$ for some $c\in\R$, or $f=\mathrm{id}$.
    \item $f=g+h$ or $f=g\cdot h$ for two (already constructed) polynomials $g$ and $h$. 
\end{enumerate}
Here $+$ and $\cdot$ denote pointwise addition and multiplication of real 
functions. We denote the set of polynomials again by $\R[x]$. Now arithmetic 
operations in $\R[x]$ are trivial, they are inbuilt in the definition. The above 
expressions take form
$${\textstyle
p(x)=\sum_{j=0}^n a_jx^j=\sum_{j=0}^n f_{a_j}\cdot\mathrm{id}^j
}
$$ 
and serve not as definitions but as canonical forms of polynomials. The abbreviation EF refers to  so called 
elementary functions in real analysis. $\R[x]$ forms (via the EF-definition) their natural and important subclass.

But we are more interested in algebraic definitions of polynomials. In Algebra polynomials in $\R[x]$ are often defined as {\em formal} expressions
$$
{\textstyle
p=p(x)=\sum_{j=0}^n a_jx^j
}
$$
with $n$ and $a_j$ as above but with $x$ being a~formal 
variable. We call it the A-{\em definition} (for example \cite[Chapter~3.1]{macl_birk} or \cite[Chapter 3.14]{vdWa}). It suffers from the same 
problem as the MA-definition and from an additional drawback: ``formal expressions'' 
are vague, what sets are they?

We conclude our review of polynomials in this section with another 
algebraic definition, the IM-{\em definition}; later we give other definitions. IM stands 
for ``infinite maps''. In it $\R[x]$ consists of all maps $p\cc\N_0=\{0,1,\ds\}\to\R$ such that 
$p(n)\ne0$ for only finitely many $n$ (for example \cite[p.~97]{lang}). The 
MA-definition problem persists and the additional drawback is that each such $p$ 
is an infinite set, but formal polynomials should be finite objects. We fix it in 
Section~\ref{sec_standPoly} by keeping only nonzero values of $p$.

In the first part of our article in Sections~\ref{sec_defini}--\ref{sec_polynomials} we develop a~new formal theory of multivariate 
polynomials over commutative rings. It is a~formal version of the EF-definition. In
Section~\ref{sec_defini} we summarize notation and notions used later. In 
Section~\ref{sec_standPoly} we review in detail the standard FM-{\em definition} (``finite maps'') of 
polynomials. We extend it to multivariate polynomials over general rings. We include this ``standard'' 
material because we did not find any sufficiently rigorous and detailed reference for it and 
because we need standard polynomials later for  the map $\Psi$ in 
Section~\ref{sec_polynomials}. In Section~\ref{sec_terms} we review terms 
with arbitrary arity function and evaluation maps on algebras corresponding to terms. 

In Section~\ref{sec_polynomials} we introduce formal polynomials, with $n$ 
variables and with arbitrarily many variables, over commutative rings. They 
are s-equivalence classes of ring terms. These terms are built 
from variables, from two symbols for addition and multiplication and from constant symbols representing elements of the ring. Two terms 
are s-equivalent if they yield the same standard polynomial when one substitutes ``true'' variables for the formal variables and elements of the ring for the constant symbols. Two terms are 
f-equivalent if one can be transformed in the other by a~series of elementary 
transformations. Each of these transformations of terms (see Definition~\ref{def_elemEquiv}) replaces a~subterm in a~term with another 
subterm. For example, a~subterm $(u\cdot v)$, where $u$ and $v$ are terms, is 
replaced with $(v\cdot u)$. A~main result of our article is 
Theorem~\ref{thm_sEqIsfEq} that s-equivalence and f-equivalence coincide.

The second part of our article in Sections~\ref{sec_BooTermPoly}--\ref{sec_PM} is 
devoted to Boolean Independence Theorem~\ref{thm_indepThm}. In 
Section~\ref{sec_BooTermPoly} we review Boolean algebras and 
Boolean formulas, which we call (Boolean) terms, and introduce evaluation maps. In 
Section~\ref{sec_DNFterms} we again introduce f-equivalence of 
terms and review DNF (disjunctive normal form) terms. Section~\ref{sec_BoolInde} contains 
Definition~\ref{def_probSpace} of finitary probability spaces, 
Definition~\ref{def_Indep} of mutual independence of two tuples of events 
and another main result of our article, 
Theorem~\ref{thm_indepThm}. It says that two events which are values 
of the evaluation map on two terms and two mutually independent tuples of 
events are independent. In Section~\ref{sec_PM} we show in an example how this 
theorem justifies ``obvious'' independence of events in an application of 
the Lov\'asz Local Lemma to colorings of hypergraphs. In the last 
Section~\ref{sec_conclCom} we give 
psychologising concluding remarks.

\section{Notation and notions}\label{sec_defini}

We use standard set-theoretic notation, for example $\sus$ denotes non-strict 
inclusion of sets and $\mathcal{P}(y)=\{x\;|\;x\sus y\}$ is the power set of 
the set $y$. $\N=\{1,2,\ds\}$ is the set of natural numbers, $\N_0=\N\cup\{0\}=\{0,1,\ds\}$ 
and $\Z$ is the ring of integers. For $n\in\N_0$ we set $[n]=\{1,2,\ds,n\}$, with $[0]=\emptyset$. For a~finite set $X$ we denote by $|X|\in\N_0$ its 
cardinality.

Let $X$ and $Y$ be sets. Notation $f\cc X\to Y$ indicates that $f\sus X\times Y$ and 
that for every $x\in X$ there is a~unique $y\in Y$ with $(x,y)\in f$, written 
$f(x)=y$. We use also the alternative notation $f[x]=y$ because round brackets are needed elsewhere. We say that $f$ 
is a~function, or a~map, from $X$ to $Y$.

For a~nonempty set $A$, the set
$$
A^*=\{u\;|\;u\cc[n]\to A,\,n\in\N_0\}
$$
comprises all finite sequences with terms in $A$. The elements of $A^*$ are called {\em words (over the 
alphabet $A$)} and we write them as
$$
u=u_1\,u_2\,\ds\,u_n\text{ where the $u_i=u(i)$ are the {\em letters} of $u$}\;.
$$
The {\em length} $|u|\in\N_0$ of the word $u$ is the number $n$. 
For $n=0$ we have the unique empty word $\emptyset$. For $n=1$ we assume that 
it is known from the context
if the notation $u_1$ means the singleton word $\{(1,u(1))\}$ of length $1$, or 
just the set $u(1)$. Let $k,n\in\N$. A~{\em subword $v$} of $u$ is any word $v=v_1v_2\ds v_k$ such that 
$$
v_1=u_i,\, v_2=u_{i+1},\, \ds,\, v_k=u_{i+k-1} 
$$
for some index $i$ with $1\le i\le n$ and $i+k-1\le n$. It also means the 
subsequence (subset) $u_iu_{i+1}\ds u_{i+k-1}$ of $u$. The empty word is 
a~subword of any word.

The {\em concatenation $w=uv$} of two words $u=u_1\ds u_k$ and $v=v_1\ds v_l$ over $A$ ($k,l\in\N_0$) is the word $w=w_1\ds w_{k+l}$ over $A$ such that $w_i=u_i$ for $i\in[k]$ and $w_i=v_{i-k}$ for
$i\in[k+l]\setminus[k]$. This is an associative (in general not 
commutative) operation on $A^*$. Thus for $k\in\N$ and any $k$-tuple 
$u_{(i)}\in A^*$, $i\in[k]$, of words, the {\em concatenated word}
$$
u_{(1)}\,u_{(2)}\,\ds\,u_{(k)}\in A^*
$$
is unambiguously defined.

By a~{\em ring} 
$$
R=(R,\,\oplus,\,\odot,\,0_R,\,1_R)
$$ 
we mean a~commutative ring with the identity $1_R$. It is a~quintuple 
of a~set  $R$, two binary operations
$$
\oplus,\,\odot \cc R\times R\to R
$$
on it and two distinct constants $0_R,1_R\in R$, so called {\em neutral elements}. It satisfies the following well known ring axioms. The operations $\oplus$ and $\odot$ are associative
and commutative, for every $r\in R$ one has that $0_A\oplus r=r$ and $1_A\odot r=r$, every element $r\in R$ has an (additive) inverse $-r\in R$ satisfying 
$r\oplus(-r)=0_R$, and $\odot$ is distributive to $\oplus$, which means that for every $r,r',r''\in R$ the equality 
$$
r\odot(r'\oplus r'')=(r\odot r')\oplus(r\odot 
r'') 
$$
holds. We denote ring operations by $\oplus$ and $\odot$ because $+$ and 
$\cdot$ are used elsewhere.

If $R$ is a~ring, an {\em $R$-algebra} is any algebra 
$$
X=(X,\,\{h_j\;|\;j\in\{+,\,\cdot,\,0,\,1\}\cup R\})
$$
with the base set $X\ne\emptyset$ and operations $h_j$ on it such that $h_+$ and $h_{\cdot}$ are binary operations, $h_0,h_1\in X$ and $h_r\in X$ for $r\in R$ are constants such that $h_0=h_{0_R}$ and $h_1=h_{1_R}$, and the subalgebra
$$
(\{h_r\;|\;r\in R\},\,h_+,\,h_{\cdot},\,h_0,\,h_1)
$$
is a~ring isomorphic to $R$. For example, any ring $R$ is the $R$-algebra
$$
(R,\{\oplus,\,\odot,\,0_R,\,1_R\}\cup\{h_r\;|\;r\in R\})
$$
with the constant $h_r$ realized by the element $r$. 

\section[Standard polynomials over commutative rings]{Standard polynomials}\label{sec_standPoly}

In this section we present the FM-{\em definition} (``finite maps'') of multivariate polynomials over rings. We define arithmetic 
operations with them  and show that they form rings. We explained earlier 
why we present here in detail such ``standard'' material. 

Let $n\in\N$. The set
$$
\N_0^n=\N_0\times\N_0\times\ds\times\N_0=\{\overline{m}=(m_1,\,\ds,\,m_n)\;|\;m_i\in\N_0\}
$$
comprises $n$-tuples of numbers in $\N_0$. The set
$$
\N_0^{\infty}
:=\{\emptyset\}\cup\N\cup(\N_0\times\N)\cup
(\N_0\times\N_0\times\N)\cup\ds
$$
comprises words over $\N_0$ with nonzero last letters. We define standard polynomials with $n$ variables and with arbitrarily many variables.

\begin{defi}\label{def_stanPoly}
Let $n\in\N$ and $R$ be a~ring. A~standard polynomial $p$ over $R$ with $n$ variables is a~map
$$
p\cc X\to R\setminus\{0_R\}
$$
from a~finite set $X\sus\N_0^n$ to nonzero elements of the ring. The set of all such polynomials is denoted by $R[x_1,\ds,x_n]$.
\end{defi}

\begin{defi}\label{def_stanPolyInf}
Let $R$ be a~ring. A~standard polynomial $p$ over $R$ with arbitrarily many variables is a~map
$$
p\cc X\to R\setminus\{0_R\}
$$
from a~finite set $X\sus\N_0^{\infty}$ to nonzero elements of the ring. The set of all such polynomials is denoted by $R[x_1,x_2,\ds]$.
\end{defi}
By keeping only nonzero coefficients we ensure uniqueness of these representations 
of polynomials. The equality sign $\stackrel{s}{=}$ indicates the informal A-definition style notation for 
standard polynomials. For example, with $R=\Z$ we have in $\Z[x_1,\ds,x_4]$ the standard polynomial
$$
p=\{((0,0,0,0),\,3),\;((0,2,0,0),\,-1),\;((1,0,3,0),\,1)\}\stackrel{s}{=}3-x_2^2+x_1x_3^3\;.
$$
In $\Z[x_1,x_2,\ds]$ it takes the form
$$
p=\{(\emptyset,\,3),\;((0,2),\,-1),\;((1,0,3),\,1)\}\stackrel{s}{=}3+x_3^3x_1-x_2^2\;.
$$
If $n\in\N$, $r\in R\setminus\{0_R\}$ and $\overline{m}=(m_1,\ds,m_n)\in\N_0^n$, we write
$$
rx^{\overline{m}}\stackrel{s}{=}rx_1^{m_1}x_2^{m_2}\ds x_n^{m_n}
$$
for the standard monomial $\{((m_1,\,m_2,\,\ds,\,m_n),r)\}$ in $R[x_1,\ds,x_n]$, and similarly in $R[x_1,x_2,\ds]$.

For $\overline{k},\overline{l}\in\N_0^n$ with $n\in\N$ we set
$$
\overline{k}+\overline{l}:=(k_1+l_1,\,\ds,\,k_n+l_n)\;.
$$
For $\overline{k},\overline{l}\in\N_0^{\infty}$ with $\overline{k}=(k_1,\ds,k_m)$ and $\overline{l}=
(l_1,\ds,l_n)$, $m,n\in\N_0$ and $m\le n$, we set $\overline{k}+\overline{l}:=\overline{l}$ if $m=0$ (then 
$\overline{k}=\emptyset$), and 
$$
\overline{k}+\overline{l}:=
(k_1+l_1,\,\ds,\,k_m+l_m,\,l_{m+1},\,\ds,\,l_n)
$$
if $m>0$. We define arithmetic operations on polynomials.

\begin{defi}\label{def_opePol1}
Let $R$ be a~ring, $n\in\N$ and 
$$
p\cc X\to R\setminus\{0_R\}\,\text{ and }\, q\cc Y\to R\setminus\{0_R\}
$$ be two polynomials in $R[x_1,\ds,x_n]$. Their sum  is the polynomial
$$
p\oplus_s q\cc Z\to R\setminus\{0_R\},\ Z\sus X\cup Y\;,
$$
where we put in $Z$ every $z\in X\cup Y$ such that either {\em (i)} $z$ lies only in $X$ or {\em (ii)} $z$ lies 
only in $Y$ or {\em (iii)} $z\in X\cap Y$ and $\nu:=p(z)\oplus q(z)\ne0_R$. In the case {\em (i)} we 
set $(p\oplus_s q)(z):=p(z)$, in the case {\em (ii)} $(p\oplus_s q)(z):=q(z)$ and in the case {\em (iii)} 
$(p\oplus_s q)(z):=\nu$. The product of $p$ and $q$ is the polynomial 
$$
p\odot_s q\cc T\to R\setminus\{0_R\},\ T\sus\N_0^n\;,
$$
defined by computing for every $\overline{m}\in\N_0^n$ the convolution
$$
\nu(\overline{m}):=\bigoplus_{\overline{k}+\overline{l}=\overline{m}}p(\overline{k})\odot q(\overline{l})\in R\;,
$$
where the addition $\oplus$ is in $R$ and $p$ and $q$ are extended outside their definition domains by $0_R$, and by setting 
$$
T:=\{\overline{m}\in\N_0^n\;|\;\nu(\overline{m})\ne
0_R\}
$$ 
and $(p\odot_s q)(\overline{m}):=\nu(\overline{m})$ for every $\overline{m}\in T$. Finally, we define the neutral elements by
$$
0_{R[x_1,\ds,x_n]}=\emptyset\,\text{ and }\,1_{R[x_1,\ds,x_n]}=\{((0,\,\ds,\,0),\,1_R)\}\;.
$$
\end{defi}
It is easy to see that the set $T$ is finite.

\begin{defi}\label{def_opePol2}
Let $R$ be a~ring and let 
$$
p\cc X\to R\setminus\{0_R\}\,\text{ and }\, q\cc Y\to R\setminus\{0_R\}
$$ 
be two polynomials in $R[x_1,x_2,\ds]$. The sum $p\oplus_s q$ and the product $p\odot_s q$ are defined as in 
Definition~\ref{def_opePol1}, with the only change that $\overline{m}\in\N_0^{\infty}$ and that in the product one 
uses addition of tuples of indices in $\N_0^{\infty}$. The neutral elements are 
$$
0_{R[x_1,x_2,\ds]}=\emptyset\,\text{ and }\,1_{R[x_1,x_2,\ds]}=\{(\emptyset,\,1_R)\}\;.
$$
\end{defi}
In $\oplus_s$ and $\odot_s$ the lower index $s$ indicates standard polynomials. In 
Section~\ref{sec_polynomials} we introduce addition and multiplication of formal 
polynomials. 

We show that standard polynomials form rings.

\begin{thm}\label{thm_polyRingStand}
Let $n\in\N$, $R$ be a~ring and $\mathcal{R}=R[x_1,\ds,x_n]$, resp. $\mathcal{S}=R[x_1,x_2,\ds]$. Then
$$
(\mathcal{R},\,\oplus_s,\,\odot_s,\,0_{\mathcal{R}},\,1_{\mathcal{R}}),\,\text{ resp. }\,(\mathcal{S},\,\oplus_s,\,\odot_s,\,0_{\mathcal{S}},\,1_{\mathcal{S}})\;,
$$ 
is a~ring.
\end{thm}
\proof
We discuss in detail only the ring $\mathcal{R}=R[x_1,\ds,x_n]$, $n\in\N$. It follows from commutativity and 
associativity of $\oplus$ in $R$ that such is also $\oplus_s$ in $\mathcal{R}$. It is also clear that 
$0_{\mathcal{R}}=\emptyset$ is neutral to $\oplus_s$. 
Associativity and commutativity of $\odot_s$ in $\mathcal{R}$ follow from associativity and commutativity of $+$ in $\N_0^n$ and $\odot$ in $R$. Clearly,
$$
1_{\mathcal{R}}=\{((0,\ds,\,0),\,1_R)\}
$$
is neutral to $\odot_s$. If $p\cc X\to R\setminus\{0_R\}$ is a~polynomial in 
$\mathcal{R}$, its additive inverse
$-p\cc X\to R\setminus\{0_R\}$ is given by
$$
(-p)(\overline{k}):=-p(\overline{k})\in R\setminus\{0_R\},\ \overline{k}\in X\;.
$$
The satisfaction of the distributive law is least obvious. Let $p$, $q$ and $r$ be three polynomials in $\mathcal{R}$ 
and let $\overline{k}\in\N_0^n$. We extend polynomials outside their definition domains by $0_R$. We see that indeed
\begin{eqnarray*}
\big(p\odot_s(q\oplus_s r)\big)(\overline{k})&=&\bigoplus_{\overline{l}+\overline{m}=\overline{k}}p(\overline{l})\odot(q\oplus_s r)(\overline{m})\\
&=&\bigoplus_{\overline{l}+\overline{m}=\overline{k}}
p(\overline{l})\odot(q(\overline{m})\oplus r(\overline{m}))\\
&=&\bigoplus_{\overline{l}+\overline{m}=\overline{k}}\big((p(\overline{l})\odot q(\overline{m}))\oplus(p(\overline{l})\odot r(\overline{m}))\big)\\
&=&\bigoplus_{\overline{l}+\overline{m}=\overline{k}}p(\overline{l})\odot q(\overline{m})\oplus
\bigoplus_{\overline{l}+\overline{m}=\overline{k}}p(\overline{l})\odot r(\overline{m})\\
&=&(p\odot_s q)(\overline{k})\oplus(p\odot_s r)(\overline{k})=\big((p\odot_s q)\oplus_s(p\odot_s r)\big)(\overline{k})\;.
\end{eqnarray*}
The interested reader easily justifies every step in the computation.

Similarly one proves that $\mathcal{S}=R[x_1,x_2,\ds]$ is a~ring. The only changes are that
$1_{\mathcal{R}}=\{(\emptyset,1_R)\}$,
that $\overline{k}\in\N_0^{\infty}$ and that addition of tuples of indices takes place in $\N_0^{\infty}$.
\eproof

In the next definition we show that the rings $R[x_1,\ds,x_n]$ and 
$R[x_1,x_2\ds]$ are $R$-algebras.

\begin{defi}\label{def_stanPoly3}
Let $n\in\N_0$ and $R$ be a~ring. For $n=0$ the ring of standard polynomials over $R$ with $n$ (i.e., without) variables is the $R$-algebra $R$, as explained earlier. For $n\ge1$ it is the $R$-algebra
$$
\big(R[x_1,\,\ds,\,x_n],\,\{\oplus_s,\,\odot_s,\,0_{R[x_1,\ds,x_n]},\,1_{R[x_1,\ds,x_n]}\}\cup\{c_r\;|\;r\in R\}\big)
$$
with the constant $c_r$ realized by $\{((0,0,\ds,0),r)\}$ for $r\ne0_R$ and by  $\emptyset$ for $r=0_R$.
The ring of standard polynomials over $R$ with arbitrarily many variables is the $R$-algebra
$$
\big(R[x_1,\,x_2,\,\ds],\,\{\oplus_s,\,\odot_s,\,0_{R[x_1,x_2\ds]},\,1_{R[x_1,x_2,\ds]}\}\cup\{c_r\;|\;r\in R\}\big)
$$
with the constant $c_r$ realized by $\{(\emptyset,r)\}$ for $r\ne0_R$ and by $\emptyset$ for $r=0_R$.
\end{defi}

\section{Terms and evaluation maps}\label{sec_terms}

For the definition of formal polynomials in the next section we need terms and 
evaluation maps and we introduce them now. We define terms completely generally, even 
if later we only use functions with arity at most $2$.

Let 
$$
V=\{x_i\;|\;i\in\N\}=\{x_1,\,x_2,\,\ds\}
$$ 
be a~countable set of distinct {\em variables} and let $F\ne\emptyset$ be a~set of {\em functional symbols}. We assume that $F\cap V=\emptyset$. We also assume that the comma , and the right and left 
bracket $)$ and $($ are not elements of $F\cup V$. An {\em arity (function)} is any map
$$
a\cc F\to\N_0\;.
$$ 
If $a[c]=0$, we say that $c\in F$ is a~{\em constant}. We define terms, subterms, constant terms and atomic terms.

\begin{defi}\label{def_terms}
Let $V$, $F$ and $a$ be as above. The $a$-terms are exactly the nonempty words over the alphabet 
$$
A_F:=V\cup F\cup\{\,,\ )\ (\,\}
$$ 
that can be obtained by finitely many applications of four rules. 
\begin{enumerate}
    \item Every constant $c\in F$ and every variable $x\in V$ is a~term.
    \item If $f\in F$ has arity $a[f]=1$ and $t$ is a~term then the concatenated word $f(t)$ is a~term.
    \item If $f\in F$ has arity $a[f]=2$ and $t$ and $u$ are terms then the concatenated word $(tfu)$ is a~term.
    \item If $f\in F$ has arity $a[f]=k\ge3$ and $t_1$, $t_2$, $\ds$, $t_k$ are $k$ terms then 
    the concatenated word $f(t_1,t_2,\ds t_k)$ is a~term.
\end{enumerate}
An $a$-term $u$ is a~subterm of an $a$-term $t$ if $u$ is a~subword of the word $t$. 
A~constant term contains no variable. A~term is atomic if it is 
a~variable or a~constant 
\end{defi}
We remark several things. The type of terms is determined by the 
arity function $a$ and therefore we speak of $a$-terms. 
We could use Polish notation and do without the comma and brackets, but we prefer notation close to practice. In the definition of the alphabet $A_F$ we do not use for clarity the comma as a~separator 
and separate the three elements by two spaces. But elsewhere we freely use the comma as a~separator and 
trust reader's ability to distinguish two modes of its use, either as an element of $A_F$ or as a~meta-separator.
The same applies to the brackets $)$ and $($, and to the variables $x_i$.
The respective lengths of the concatenated words in rules 1--4 of Definition~\ref{def_terms} are 
$$
1,\ |t|+3,\ |t|+|u|+3\,\text{ and }\,|t_1|+|t_2|+\ds+|t_k|+k+2\;. 
$$
In the concatenated words of rules 1--4 the symbols $c$, $x$, $f$, $t$, $u$, $t_1$, $t_2$, `$\ds$' and $t_k$
are syntactic metavariables, but the symbols $($, $)$ and $,$ are direct mentions of elements
of $A_F$. In rule~4 the symbol $\ds$ of elipsis appears twice in two different meanings. 
The latter one indicates a~certain word over $A_F$ with length $\sum_{i=3}^{k-1}|t_i|+k-3$; for $k=3$ 
this word is empty. Thus the missing comma before $t_k$ is not a~typo but a~feature. Finally, instead of speaking of `finitely many applications $\ds$' 
we could state Definition~\ref{def_terms} more rigorously by defining $a$-terms as the final letters in certain `generating words' 
$$
U\in(A_F^*)^*\;,
$$
but the present level of rigor is sufficient. 

We denote the set of all $a$-terms by $T_a$. For $W\sus V$ we denote by $T_a[W]$ the set of terms $t\in T_a$ using only the variables in $W$.
Thus if $t\in T_a[W]$ and $x\in V\setminus W$ then $x$ does not appear in $t$, and if $x\in W$ then $x$ may 
or may not appear in $t$. For $W=\{x_1,\ds,x_n\}$ we write $T_a[x_1,\ds,x_n]$. For $W=\emptyset$ the set $T_a[\emptyset]$ consists of constant $a$-terms. 

Suppose that $F$ is a~set of functional symbols with the arity $a$ and that 
$$
\mathcal{A}=(X,\,\{h_r\;|\;r\in F\})
$$ 
is an algebra with the signature $a$. Thus for any $r\in F$,
$$
h_r\cc X^{a(r)}=X\times X\times\ds\times X\to X
$$
is an $a(r)$-ary operation on the base set  $X$. Let $n\in\N$. We associate with every $a$-term $t\in T_a[x_1,\ds,x_n]$ the {\em evaluation map}
$$
\Phi_{\mathcal{A}}(t;\,y_1,\,\ds,\,y_n)\cc X^n=X\times X\times\ds\times X\to X\;.
$$
Informally, for a~given $n$-tuple $\overline{b}=(b_1,\ds,b_n)\in X^n$ we substitute for every occurrence of $x_i$ in $t$ the element 
$b_i\in X$, perform all operations $h_r$ indicated in $t$ by the functional symbols $r$, and output the value
$$
\Phi_{\mathcal{A}}(t;\,\overline{b})=\Phi_{\mathcal{A}}(t;\,b_1,\,\ds,\,b_n)\in X\;.
$$
For the precise definition of $\Phi_{\mathcal{A}}$ we need the following {\em unique reading lemma} for terms. Its proof is left to the interested reader. 

\begin{thm}\label{thm_uniqRead}
Let $V$, $F$ and $a$ be as above and let $t\in T_a$. Then exactly one of the following four cases occurs.
\begin{enumerate}
    \item There is a~unique constant $c\in F$ or a~unique variable $x\in V$ such that $t=c$ or $t=x$.
    \item There is a~unique functional symbol $f\in F$ with $a[f]=1$ and a~unique term 
    $u\in T_a$ such that $t$ is the concatenated word $f(u)$. 
    \item There is a~unique functional symbol $f\in F$ with $a[f]=2$ and two unique terms 
    $u,v\in T_a$ such that $t$ is the concatenated word $(ufv)$. 
    \item There is a~unique number $k\in\N$ with $k\ge3$, a~unique functional symbol $f\in F$ 
    with $a[f]=k$ and $k$ unique terms $u_i\in T_a$, $i\in[k]$, such that $t$ is the concatenated 
    word $f(u_1,u_2,\ds u_k)$. The elipsis $\ds$ indicates a~word over $A_F$ of length 
    $k-3+\sum_{i=3}^{k-1}|u_i|$.
\end{enumerate}
If $t'$ is a~subterm of $t$, $t'\ne t$, then $t'$ is either a~subterm of $u$ in case~2 or of $u$ in 
case~3 or of $v$ in case~3 or of a~single $u_i$ in case~4.
\end{thm}
In mathematical logic, unique reading lemmas --- for terms and formulas --- are basic syntactic results. 

It is is easy to prove by induction on the structure of terms the next result on 
substitution; we omit the proof.

\begin{cor}\label{cor_subsTerm}
Any replacement of a~subterm in an $a$-term by an $a$-term produces an $a$-term.
\end{cor}
With the help of the last theorem we rigorously define evaluation maps.

\begin{defi}\label{def_evalMap}
Suppose that $n\in\N_0$, that
$$
\mathcal{A}=(X,\,\{h_r\;|\;r\in F\})
$$ 
is an algebra with the signature $a\cc F\to\N_0$ and that
$t\in T_a[x_1,\ds,x_n]$ (if $n=0$, $t$ is a~constant term). We define 
by induction on the structure of $t$ described in Theorem~\ref{thm_uniqRead} the evaluation map
$$
\Phi_{\mathcal{A}}[t;\,y_1,\,\ds,\,y_n]\cc X^n\to X
$$
(for $n=0$ we have $\Phi_{\mathcal{A}}[t]\in X$).
For $\overline{b}=(b_1,\ds,b_n)\in X^n$ we define the value
$$
\nu=\Phi_{\mathcal{A}}[t;\,\overline{b}]=\Phi_{\mathcal{A}}[t;\,b_1,\,\ds,\,b_n]\in X
$$
as follows. If $t=c\in F$ is a~constant we set $\nu:=h_c\in X$. If $t=x_i\in V$, 
$i\in[n]$, we set $\nu:=b_i\in X$. In the unary case when $t=f(u)$ we set by induction
$$
\nu:=h_f(\Phi_{\mathcal{A}}[u;\,\overline{b}])\in X\;.
$$
In the binary case when $t=(ufv)$ we set by induction
$$
\nu:=h_f(\Phi_{\mathcal{A}}[u;\,\overline{b}],\,\Phi_{\mathcal{A}}[v;\,\overline{b}])\in X\;.
$$
Finally, in the $k$-ary case, $k\ge3$, when $t=f(u_1,u_2,\ds u_k)$ we set by induction
$$
\nu:=h_f(\Phi_{\mathcal{A}}[u_1;\,\overline{b}],\,\Phi_{\mathcal{A}}[u_2;\,\overline{b}],\,\ds,\,\Phi_{\mathcal{A}}[u_k;\,\overline{b}])\in X\;.
$$
\end{defi}
By Theorem~\ref{thm_uniqRead} this is a~correct definition, always exactly one of the described cases holds of $t$.

For later reductions of terms to standard forms we need the {\em depth $\mathrm{dep}[t]\in\N_0$} of a~term $t$.

\begin{defi}\label{def_depth}
We define the map $\mathrm{dep}\cc T_a\to\N_0$ by induction on the structure of terms (based on Theorem~\ref{thm_uniqRead}) as follows. If $t$ is atomic, we set $\mathrm{dep}[t]:=0$. In the $k$-ary case with $k\ge1$ when $t$ is formed by the terms $t_1$, $\ds$, $t_k$ we set
$$
\mathrm{dep}[t]:=1+\max(\{\mathrm{dep}[t_1],\,\ds,\,\mathrm{dep}[t_k]\})\;.
$$
\end{defi}
For example, the term $(x_2+(c\cdot x_1))$, where $c$ is a~constant and $+$ 
and $\cdot$ are binary functions, has depth $2$.

\section[s-equivalence and 
f-equivalence of ring terms]{s-equivalence and f-equivalence of ring terms}\label{sec_polynomials}

In this section we give the F-{\em definition} (``formal'') of polynomials. We introduce s-equivalence and 
f-equivalence of ring terms. Formal polynomials over 
rings are s-equivalence classes of ring terms. Our main result is 
Theorem~\ref{thm_sEqIsfEq} which says that both equivalences of ring terms 
coincide.

Let $R=(R,\oplus,\odot,0_R,1_R)$ be a~ring. We associate with it the functional symbols $F_R:=\{+,\cdot,0,1\}\cup\{c_r\;|\;r\in R\}$ and the arity $a_R\cc F_R\to\N_0$,
$$
a_R[+]=a_R[\cdot]=2\,\text{ and }\,a_R[0]=a_R[1]=a_R[c_r]=0,\ r\in R\;.
$$
The constants $c_r$ formally represent elements of $R$. We fix notation and terminology for $a_R$-terms.

\begin{defi}\label{def_aRterms}
Let $R$ be a~ring and $a_R$ be the above arity function. In this section (ring) terms mean $a_R$-terms. We denote their sets by $$
T:=T_{a_R}\,\text{ and }\,T[x_1,\ds,x_n]:=T_{a_R}[x_1,\ds,x_n]\;.
$$
Instead of terms we could speak of protopolynomials, it is fitting but too long.
\end{defi}

We define two binary and two constant operations on $T$. 

\begin{defi}\label{def_aritFormPoly}  
Let $t,u\in T$ be terms. Their sum $t\oplus_f u$ and product $t\odot_f u$ are the concatenated words and terms
$$
t\oplus_f u:=(t+u)\in T\,\text{ and }\,t\odot_f u:=(t\cdot u)\in T\;.
$$
The neutral elements are the terms $0_T:=0\in T$ and $1_T:=1\in T$.
\end{defi}
The lower index ${}_f$ refers to formal nature of this arithmetic. Note that every set of terms $T[x_1,\ds,x_n]$, $n\in\N_0$, is closed 
to all four operations.

Terms are formal polynomials over $R$ in the most general form, but 
they do not form rings. To get rings we factorize $T$
by s-equivalence. We call two terms s-equivalent if substitutions, realized by the map $\Psi\cc T\to R[x_1,x_2,\ds]$ below, of true variables for formal ones and of elements of $R$ for ring constants yield the same standard polynomial.

\begin{defi}\label{def_psi}
Let $t$ be any term and the number $n\in\N_0$ be minimum such that $t\in T[x_1,\ds,x_n]$. We set
$$
\Psi[t]:=\Phi_{\mathcal{R}}[t;\,\overline{v}]\in R[x_1,\,x_2,\,\ds]
$$
where $\mathcal{R}=R[x_1,x_2,\ds]$ is the $R$-algebra of standard polynomials over $R$ with arbitrarily many variables (see Definitions~\ref{def_evalMap} 
and \ref{def_stanPoly3} and Theorem~\ref{thm_polyRingStand}) and $\overline{v}=(v_1,\ds,v_n)$ is the $n$-tuple of standard polynomials
$$
v_i:=\{(\underbrace{(0,\,\ds,\,0,\,1)}_{\text{$i$-tuple}},\,1_R)\}\stackrel{s}{=}1_R x_i\stackrel{s}{=}x_i\in R[x_1,\,x_2,\,\ds],\ i\in[n]\;.
$$
For $n=0$, when $t$ is a~constant term, 
$$
\Psi[t]=\{(\emptyset,\,\Phi_R[t])\}\stackrel{s}{=}\Phi_R[t]\,\text{ if }\,\Phi_R[t]\in R\setminus\{0_R\}
$$
and $\Psi[t]=\emptyset\stackrel{s}{=}0_R$ if $\Phi_R[t]=0_R$.
\end{defi}
For 
example, on the terms
$$
t:=(x_1+x_2),\ u:=(x_2\cdot(x_1+x_5))\,\text{ and }\,v:=((1+c_{1_R})\cdot c_{1_R})
$$
the map $\Psi$ has the values
\begin{eqnarray*}
\Psi[t]&=&\{((1),\,1_R),\,((0,1),\,1_R)\}\stackrel{s}{=}x_2+x_1\;,\\
\Psi[u]&=&\{((0,1,0,0,1),\,1_R),\,((1,1),\,1_R)\}\stackrel{s}{=}x_1x_2+x_5x_2\;\text{ and}\\
\Psi[v]&=&\{(\emptyset,\,2_R)\}\stackrel{s}{=}2_R,\text{ assuming that $2_R\ne0_R$}\;.
\end{eqnarray*}

\begin{defi}\label{def_sEquiv}
We say that terms $t$ and $u$ are s-equivalent, written $t\stackrel{s}{\approx}u$, if $\Psi[t]=\Psi[u]$.
\end{defi}

$R[x_1,x_2,\ds]$ is generated over $R$ by the (true) variables $x_1$, $x_2$, $\ds$ and we get the following result.

\begin{prop}\label{prop_onSequiv}
Let $n\in\N_0$, $\mathcal{R}:=R[x_1,x_2,\ds]$ and $t$ and $u$ be s-equivalent terms in $T[x_1,\ds,x_n]$. Then
$$
\Phi_{\mathcal{R}}[t;\,\overline{b}]=\Phi_{\mathcal{R}}[u;\,\overline{b}]
$$
for every $\overline{b}\in R[x_1,x_2,\ds]^n$.
\end{prop}
For any $R$-algebra $X$, number $n\in\N_0$ and term $t\in 
T[x_1,\ds,x_n]$ we have the evaluation map
$$
\Phi_X[t;\,\overline{y}]=
\Phi_X[t;\,y_1,\,\ds,\,y_n]\cc X^n\to X\;.
$$
It is well known that we may have $\Phi_X[t;\overline{y}]=\Phi_X[u;\overline{y}]$ for $t\not\stackrel{s}{\approx}u$. One of the simplest examples comes for $R=X=\Z_2$ and $n=1$ from
$$
t=((x_1\cdot x_1)+x_1)\,\text{ and }\,u=0\;.
$$
Then $\Psi[t]\stackrel{s}{=}x_1^2+x_1$ and $\Psi[u]\stackrel{s}{=}0_R$ and 
$\Phi_X[t;y_1]=\Phi_X[u;y_1]$ is the constantly zero (two-element) map. 

The next proposition and definition give the F-definition of formal polynomials 
over rings. They are s-equivalence classes of terms.

\begin{prop}\label{prop_}
Let $n\in\N_0$ and $R=(R,\oplus,\odot,0_R,1_R)$ be a~ring. Then
$$
\big(T[x_1,\,\ds,\,x_n]/\!\stackrel{s}{\approx},\,\oplus_f,\,\odot_f,\,0_T,\,1_T\big)\,\text{ and }\,\big(T/\!\stackrel{s}{\approx},\,\oplus_f,\,\odot_f,\,0_T,\,1_T\big)
$$
are commutative rings with the identity $1_R$. The former ring is isomorphic via $\Psi$ to $R[x_1,\ds,x_n]$, and the latter to $R[x_1,x_2,\ds]$.
\end{prop}
\proof
We consider $T[x_1,\,\ds,\,x_n]/\!\stackrel{s}{\approx}$, for $T/\!\stackrel{s}{\approx}$ we argue similarly. Let $t$, 
$u$ and $v$ be terms in $T[x_1,\ds,x_n]$ with $t\stackrel{s}{\approx}u$. Since
$$
\Psi[t\oplus_f v]=\Psi[(t+u)]=\Psi[t]\oplus_s\Psi[v]=\Psi[u]\oplus_s\Psi[v]=\Psi[u\oplus_f v]\;,
$$
$t\oplus_f v\stackrel{s}{\approx}u\oplus_f v$. Similarly for 
$\odot_f$. In the same way commutativity of $\oplus_s$ and $\odot_s$ implies that 
also $\oplus_f$ and $\odot_f$ are commutative. Thus both operations 
$\oplus_f$ and $\odot_f$ respect the equivalence $\stackrel{s}{\approx}$. If $u$ is a~term then ($\mathcal{R}=R[x_1,\ds,x_n]$)
$$
\Psi[0_T\oplus_f u]=\Psi[(0+u)]=\Psi[0]\oplus_s\Psi[u]=
0_{\mathcal{R}}\oplus_s\Psi[u]=\Psi[u]\;.
$$
It follows that the s-equivalence class of $0_T$ is neutral to $\oplus_f$. 
Similarly the s-equivalence class of $1_T$ is neutral to $\odot_f$. Associativity of $\oplus_f$ and $\odot_f$ and the 
distributive law follow in the same way from the 
fact that $\oplus_s$ and $\odot_s$ have these properties. It remains to show that
every s-equivalence class of terms has an additive inverse. Let $t$ be any term. Since (again $\mathcal{R}=R[x_1,\ds,x_n]$)
\begin{eqnarray*}
\Psi[t\oplus_f(c_{-1_R}\cdot t)]&=&    
\Psi[(t+(c_{-1_R}\cdot t))]\\
&=&\Psi[t]\oplus_s(-1_{\mathcal{R}})\odot_s\Psi[t]=0_{\mathcal{R}}=\Psi[0_T]\;,
\end{eqnarray*}
this inverse is the s-equivalence class of the term $(c_{-1_R}\cdot t)$.

From the previous part of the proof we see at once that 
$$
\Psi\cc T[x_1,\,\ds,\,x_n]/\!\stackrel{s}{\approx}\,\to R[x_1,\,\ds,\,x_n]
$$
is a~ring homomorphism. It is clearly surjective and it is trivially injective: 
if $\Psi[t]=\Psi[u]$ then in $T[x_1,\ds,x_n]/\!\stackrel{s}{\approx}$ the terms $t$ and $u$ are the same. Thus $\Psi$ is a~ring isomorphism.
\eproof
\vspace{-3mm}
\begin{defi}
Let $R$ be a~ring. Formal polynomials over $R$ are the elements of the rings in 
the previous proposition, namely the s-equivalence classes in 
$T[x_1,\ds,x_n]/\!\stackrel{s}{\approx}$  and in $T[x_1,x_2,\ds]/\!\stackrel{s}{\approx}$. 
\end{defi}

Besides using the s-equivalence, we can find out if two terms (protopolynomials) 
determine the same polynomial, which is  only written in two 
different ways, by using the {\em formal equivalence}, abbreviated {\em f-equivalence}. We define it by means of {\em elementary transformations}.

\begin{defi}\label{def_elemEquiv}
Let $R=(R,\oplus,\odot,0_R,1_R)$ be a~ring. Two terms $t$ and $t'$ are 
(elementary) transforms of one another, written $t\sim t'$ and $t'\sim t$, if for 
some $r,s\in R$ and $u,u',u''\in T$ the term $t'$ arises from $t$ by 
replacing a~subterm $v$ of $t$ with the subterm $v'$ by one of the twelve rules below.
\begin{description}
    \item[{\rm ET$-1$}] $v=0$ and $v'=c_{0_R}$ 
    \item[{\rm ET0}] $v=1$ and $v'=c_{1_R}$
    \item[{\rm ET1}] $v=(c_r+c_s)$ and $v'=c_{r\oplus s}$.
    \item[{\rm ET2}] $v=(c_r\cdot c_s)$ and $v'=c_{r\odot s}$.
    \item[{\rm ET3}] $v=(0+u)$ and $v'=u$.
    \item[{\rm ET4}] $v=(1\cdot u)$ and $v'=u$
    \item[{\rm ET5}] $v=(u+u')$ and $v'=(u'+u)$.
    \item[{\rm ET6}] $v=(u\cdot
     u')$ and $v'=(u'\cdot u)$.
     \item[{\rm ET7}] $v=(u+(u'+u''))$ and $v'=((u+u')+u'')$.
     \item[{\rm ET8}] $v=(u\cdot(u'\cdot u''))$ and $v'=((u\cdot u')\cdot u'')$.
     \item[{\rm ET9}] $v=(u\cdot(u'+u''))$ and $v'=((u\cdot u')+(u\cdot u''))$.
     \item[{\rm ET10}] $v=(0\cdot u)$ and $v'=0$.
\end{description}
\end{defi}
The relation $\sim$ is obviously symmetric. Notation $t\stackrel{i}{\sim}t'$ means that 
$t\sim t'$ via $\mathrm{ET}i$. Transformations ET1--ET9 correspond to ring 
operations and ring 
axioms. ET10 corresponds to the fact that in any ring multiplication by zero yields zero. ET$-1$ and ET0 were added additionally.

We define f-equivalence of terms.

\begin{defi}\label{def_fekviv}
Let $t$ and $u$ be two terms. We say that $t$ and $u$ are f-equivalent, written $t\stackrel{f}{\approx}u$, 
if there exist $k+1$ terms $t_0$, $t_1$, $\ds$, $t_k$, $k\in\N$, such that 
$$
t=t_0\sim t_1\sim\ds\sim t_k=u\ (\text{i.e., $t_{i-1}\sim t_i$ for $i\in[k]$})\;.
$$
So one can transform $t$ in $u$ by a~series of transformations listed in Definition~\ref{def_elemEquiv}.
\end{defi}
The relation $\stackrel{f}{\approx}$ is an equivalence. For example, $((x_7+c_{1_R})\cdot x_3)\stackrel{f}{\approx}((x_3\cdot x_7)+x_3)$ because
\begin{eqnarray*}
((x_7+c_{1_R})\cdot x_3)&\stackrel{\mathrm{6}}{\sim}&
(x_3\cdot (x_7+c_{1_R}))\stackrel{\mathrm{9}}{\sim}
((x_3\cdot x_7)+(x_3\cdot c_{1_R}))\\
&\stackrel{\mathrm{6}}{\sim}&((x_3\cdot x_7)+(c_{1_R}\cdot x_3))\stackrel{\mathrm{0}}{\sim}((x_3\cdot x_7)+(1\cdot x_3))\\
&\stackrel{4}{\sim}&((x_3\cdot x_7)+x_3)\;.
\end{eqnarray*}
We show that the effect of ET10 cannot be obtained via other 
transformations.

\begin{prop}\label{prop_formMultZero}
Suppose that $t$ is a~non-constant term, that the numbers $i_1$, $\ds$, $i_j$, $j\in\N$, lie in $\{-1,\ds,10\}$ and that $t_0$, 
$\ds$, $t_j$ are terms such that     
$$
t=t_0\stackrel{i_1}{\sim}
t_1\stackrel{i_2}{\sim}\ds
\stackrel{i_j}{\sim}t_j=0\;.
$$
Then there is a~$k\in[j]$ such that $i_k=10$.
\end{prop}
\proof
This is clear from the fact that non transformation different from ET10 
erases variables.
\eproof

We show that f-equivalence of terms implies their s-equivalence.

\begin{prop}\label{prop_feqImplseqPol}
Two f-equivalent terms $t$ and $u$ are s-equivalent.
\end{prop}
\proof
It suffices to show that
if $t\stackrel{i}{\sim}u$ for some $i\in\{-1,\ds,10\}$ then $\Psi[t]=\Psi[u]$.  For example, let $i=9$. 
Then there exist a~subterm $v$ of $t$ and terms $w$, $w'$ and $w''$ such that $v=
(w\cdot(w'+w''))$ and $u$ arises from $t$ by replacing $v$ with the subterm $v':=
((w\cdot w')+(w\cdot w''))$. By Theorem~\ref{thm_uniqRead} and Definition~\ref{def_psi}, the value $\Psi[t]$, resp. $\Psi[u]$, is determined 
by the values of $\Psi$ on the proper subterms of $t$, resp. $u$, disjoint to $v$, resp. $v'$, and by the value $\Psi[v]$, resp.
$\Psi[v']$. The former corresponding values are trivially equal because those subterms in 
$t$, resp. $u$, are the same. But also $\Psi[v]=\Psi[v']$ because the 
distributive law holds in the ring $R[x_1,x_2,\ds]$. Thus $\Psi[t]=\Psi[u]$.
For $i\ne9$, other transformations, the argument is very similar, the justification of the equality $\Psi[v]=\Psi[v']$ is always easily supplied.
\eproof

Multiplication by zero plays an important role in our theory. In a~short 
intermezzo we therefore generalize the well known result that in any ring $R=
(R,\oplus,\odot,0_R,1_R)$ for every $r\in R$ one has that $r\odot 0_R=0_R$. 

If $\oplus$ is a~binary operation on a~set $X$ and $x\in X$, we call $x$ {\em left 
cancellable in $\oplus$} if the map
$$
X\ni y\mapsto x\oplus y\in X
$$ 
is injective. For example, every $n\in\N_0$ is left cancellable in the usual addition $+$ on 
$\N_0$, one proves it by induction. Or any element in a~group is left cancellable in the group operation, one proves it by inverses. An element 
$0_X\in X$ is {\em right neutral in $\oplus$} if for every $x\in X$, 
$$
x\oplus 0_X=x\;.
$$ 
Another binary operation $\odot$ on $X$ is {\em left distributive to $\oplus$} if for every $x,y,z\in X$, 
$$
x\odot(y\oplus z)=(x\odot y)\oplus(x\odot z)\;.
$$

\begin{thm}\label{thm_multByZero}
Suppose that $\oplus$ and $\odot$ are binary operations on a~set $X$, with $\odot$ left distributive to $\oplus$, and that $0_X\in X$ is right neutral in $\oplus$. Then for any $x\in X$ the element
$x\odot 0_X$ is not left cancellable in $\oplus$ or $x\odot 0_X=0_X$.
\end{thm}
\proof
Let $x\in X$ be such that $x\odot 0_X$ is left cancellable in $\oplus$. Right neutrality of $0_X$ and left distributivity of $\odot$ give that 
$$
x\odot 0_X=x\odot(0_X\oplus 0_X)=(x\odot 0_X)\oplus(x\odot 0_X)\;.
$$
Using again right neutrality of $0_X$ and canceling $x\odot 0_X$ from the left in
$$
(x\odot 0_X)\oplus 0_X=(x\odot 0_X)\oplus(x\odot 0_X)
$$
we get that $0_X=x\odot 0_X$.
\eproof

\noindent
The following corollary is immediate.

\begin{cor}\label{cor_multByzero}
Suppose that $\oplus$, $\odot$, $X$ and $0_X$ are as in the previous theorem and that every element $x\in 
X\setminus\{0_X\}$ is left cancellable in $\oplus$. Then for every $x\in X$ it holds that $x\odot 0_X=0_X$.
\end{cor}

We arrive at the first main result of our article, the theorem that s-equivalence and 
f-equivalence of terms coincide.

\begin{thm}\label{thm_sEqIsfEq}
For any ring $R$ and any $a_R$-terms $t$ and $u$,
$$
t\stackrel{s}{\approx}u\iff t\stackrel{f}{\approx}u\;.
$$
\end{thm}
Theorem~\ref{thm_sEqIsfEq} has the following practical application. It is not 
clear how to check via the elementary transformations if two terms $t$ and $u$ 
are f-equivalent (although in the proof we provide a~method). By Theorem~\ref{thm_sEqIsfEq} we can reduce it to the s-equivalence. It is not hard to see that s-equivalence of $t$ and $u$ can be checked 
(by means of Theorem~\ref{thm_uniqRead} and Definition~\ref{def_psi}) in time 
polynomial in $|t|+|u|$. 

\medskip\noindent
{\bf Proof of Theorem~\ref{thm_sEqIsfEq}. }The implication $\Leftarrow$ was established in 
Proposition~\ref{prop_feqImplseqPol}. In the rest of this section we prove the opposite implication $\Rightarrow$. We assume that $t$
and $u$ are s-equivalent terms and show that $t$ can be 
turned in $u$ by a~series of transformations ET$i$, $i\in\{-1,\ds,10\}$. To this end we first 
prove Propositions~\ref{prop_subsEquiv} and \ref{prop_subsXi} on substitution and 
$\stackrel{f}{\approx}$. Then in Proposition~\ref{prop_StandMonom} to 
Definition~\ref{def_standTerm} we introduce standard terms. In Theorem~\ref{thm_redToStand} we show how to reduce any 
term to an f-equivalent standard term. Finally, in 
Propositions~\ref{prop_redToStand1}--\ref{prop_finalStep} we prove that the 
implication $\Rightarrow$ holds for standard terms.  

First we show that f-equivalence is preserved by  substitution in subterms.

\begin{prop}\label{prop_subsEquiv}
If a~term $t'$ arises from another term $t$ by replacing a~subterm $u$ of $t$ with a~subterm $u'\stackrel{f}{\approx}u$ then $t'\stackrel{f}{\approx}t$.
\end{prop}
\proof
Let $t'$, $t$, $u$ and $u'$ be as stated. 
The word $t'$ is a~term by Corollary~\ref{cor_subsTerm}.
We may assume that $u\stackrel{i}{\sim}u'$ for some $i\in\{-1,\ds,10\}$. 
Thus $u'$ arises from $u$ by replacing a~subterm $v$ of $u$ with a~subterm $v'$ such that $v$ and $v'$ are related as 
described in ET$i$. But $v$ is also a~subterm of $t$, by the definitions of subterms and of subwords. Thus $t'$ arises from $t$ by replacing $v$ 
with $v'$ and $t\stackrel{i}{\sim}t'$. We see that $t\stackrel{f}{\approx}t'$.
\eproof

If $t$ and $u$ are terms and $x$ is a~variable,
$t[x/u]$ denotes the term obtained from $t$ by 
replacing every occurrence of $x$ in $t$ with $u$. By Corollary~\ref{cor_subsTerm}, 
$t[x/u]$ is a~term. We show that f-equivalence is preserved by such substitutions.

\begin{prop}\label{prop_subsXi}
Let $x\in V$ and let $u\stackrel{f}{\approx}u'$ 
and $t\stackrel{f}{\approx}t'$ be pairs of f-equivalent terms.  Then
$$
t[x/u]\stackrel{f}{\approx} t'[x/u']\;.
$$
\end{prop}
\proof
It suffices to show for every $i\in\{-1,\ds,10\}$ that if $t\stackrel{i}{\sim}t'$ and $u\stackrel{f}{\approx}u'$ then
$$
t[x/u]\stackrel{f}{\approx}t[x/u']\,\text{ and }\,t[x/u']\stackrel{i}{\sim}t'[x/u']\;.
$$
The former f-equivalence follows by replacing every occurrence of $x$ in $t$ with $u$ and then replacing these (disjoint) $u$-copies one by one with $u'$ by 
Proposition~\ref{prop_subsEquiv}.

We prove the latter f-equivalence. For example, let $t\stackrel{9}{\sim}t'$. Then $t'$ arises from $t$ by
replacing a~subterm $v$ of $t$ with the term $v'$, where  ($v_j\in T$) 
$$
v=(v_1\cdot(v_2+v_3))\,\text{ and }\,v'=((v_1\cdot v_2)+(v_1\cdot v_3))\;. 
$$
It follows that $t'[x/u']$ arises from $t[x/u']$ by
replacing the subterm $v[x/u']$ of $t[x/u']$ with the term $v'[x/u']$. Since 
\begin{eqnarray*}
v[x/u']&=&(v_1[x/u']\cdot(v_2[x/u']+v_3[x/u']))\,\text{ and}\\
v'[x/u']&=&((v_1[x/u']\cdot v_2[x/u'])+(v_1[x/u']\cdot v_3[x/u']))\;,
\end{eqnarray*}
$t[x/u']\stackrel{9}{\sim}t'[x/u']$. For other transformations relating $t$ and $t'$ the argument is similar. 
\eproof

Let 
$\overline{m}=(m_1,\ds,m_n)\in\N_0^{\infty}$ with $n\in\N_0$. We define a~class of terms $\mu$ called {\em $\overline{m}$-monomials}. For $n=0$ when $\overline{m}=\emptyset$ we set $\mu:=1$. For $n\ge1$, $\mu$ is any term such that 
$$
\mu\in\big(\{\,\cdot\ )\ (\,\}\cup\{x_1,\,\ds,\,x_n\}\big)^*
$$
and that for $i\in[n]$ the variable $x_i$ appears in $\mu$ exactly $m_i$ times. For example,
$$
\text{$((x_1\cdot x_3)\cdot x_3)$
and $(x_3\cdot (x_3\cdot x_1))$} 
$$
are $(1,0,2)$-monomials. We say that $\mu$ {\em has  type $\overline{m}$}. A~{\em monomial} is an $\overline{m}$-monomial. We show that $\overline{m}$-monomials are mutually f-equivalent.

\begin{prop}\label{prop_StandMonom}
Let $\overline{m}\in\N_0^{\infty}$. Any two $\overline{m}$-monomials are f-equivalent.
\end{prop}
\proof
For $\overline{m}=\emptyset$ it holds trivially and we assume that $\overline{m}\ne\emptyset$. 
We say that an $\overline{m}$-monomial $\mu$ is standard if either $\mu\in V$ or, for some $k\in\N$ with $k>1$,
$$
\mu=(x_{i_1}\cdot(x_{i_2}\cdot(\ldots\cdot(x_{i_{k-1}}\cdot x_{i_k})\ds)))
$$
and $1\le i_1\le i_2\le\ds\le i_k$. Standard $\overline{m}$-monomials are uniquely determined.  
For example, $(x_1\cdot(x_3\cdot x_3))$ is the standard 
$(1,0,2)$-monomial. We show that every $\overline{m}$-monomial $\mu$ is f-equivalent via $\stackrel{6}{\sim}$ and $\stackrel{8}{\sim}$ to a~standard monomial. Since these transformations preserve numbers of occurrences of variables, 
it is the standard $\overline{m}$-monomial and the result follows. 

We proceed by induction on the size of $\mu$. Let $i\in\N$ be the smallest index of a~variable $x_i$ appearing 
in $\mu$. If $\mu=x_i$, we are done trivially. If $\mu=
(x_i\cdot\nu)$ or $\mu=(\nu\cdot x_i)$ for a~monomial $\nu$, we are done by using 
induction for $\nu$, Proposition~\ref{prop_subsEquiv} and ET6. Else we may assume, using ET6, that $\mu=(\nu\cdot\tau)$ for monomials $\nu$
and $\tau$ such that $x_i$ appears in $\nu$ and $\nu\ne x_i$. Using induction for $\nu$ and 
Proposition~\ref{prop_subsEquiv}, we may assume that 
$\nu=(x_i\cdot\nu')$ for a~monomial $\nu'$. But then
$$
\mu=((x_i\cdot\nu')\cdot\tau)\stackrel{8}{\sim}(x_i\cdot(\nu'\cdot\tau))
$$
and we are done by using induction for $(\nu'\cdot\tau)$ and Proposition~\ref{prop_subsEquiv}.
\eproof

We determine values of the map $\Psi$ (Definition~\ref{def_psi}) on monomials.

\begin{prop}\label{prop_onMonom2}
Let $R$ be a~ring, $\overline{m}=(m_1,\ds,m_n)\in\N_0^{\infty}$, $n\in\N_0$, and let $\mu$ be an $\overline{m}$-monomial. Then $\Psi[\mu]=\{(\emptyset,1_R)\}\stackrel{s}{=}1_R$ for $n=0$ and
$$
\Psi[\mu]=\{((m_1,\,\ds,\,m_n),\,1_R)\}\stackrel{s}{=}x_1^{m_1}\ds x_n^{m_n}
$$
for $n\ge1$.
\end{prop}
\proof
For $n=0$ it holds by the definitions. It holds also for $n>0$ and $m_1+\ds+m_n=1$ 
because then $m_1=\ds=m_{n-1}=0$ and $m_n=1$, and with 
$\mathcal{R}=R[x_1,x_2,\ds]$ and 
$$
\overline{v}=\big(\{((1),1_R)\},\,\ds,\,\{((0,\,\ds,\,0,\,1),\,1_R)\}\big)\stackrel{s}{=}
(x_1,\,\ds,\,x_n)
$$ 
we get by Definition~\ref{def_psi} that
$$
\Psi[\mu]=\Psi[x_n]=\Phi_{\mathcal{R}}[x_n;\,\overline{v}]=\{((0,\,\ds,\,0,\,1),\,1_R)\}
\stackrel{s}{=}x_n\;.
$$
For $n>0$ and $m_1+\ds+m_n>1$ we proceed by induction on the sum. Let $1\le 
i_1<i_2<\ds<i_k=n$, $k\in\N$, be the indices $i_j$ such that $m_{i_j}>0$. We replace $\mu$ by the 
f-equivalent standard $\overline{m}$-monomial $(x_{i_1}\cdot\nu)$ where $\nu$ is a~standard monomial. Using Proposition~\ref{prop_feqImplseqPol} we get that $\Psi[\mu]$ equals
$$
\Psi[(x_{i_1}\cdot\nu)]=\Psi[x_{i_1}]
\odot_s\Psi[\nu]\stackrel{s}{=}x_{i_1}\odot_s x_{i_1}^{m_{i_1}-1}x_{i_2}^{m_{i_2}}\ds x_{i_k}^{m_{i_k}}\stackrel{s}{=}
x_{i_1}^{m_{i_1}}\ds x_{i_k}^{m_{i_k}}
$$
which $\stackrel{s}{=}x_1^{m_1}\ds x_n^{m_n}$.
\eproof

Let $n\in\N$. An {\em additive $n$-monomial} is any term $s$ such that 
$$
s\in\big(\{\,+\ )\ (\,\}\cup \{x_1,\,\ds,\,x_n\}\big)^*
$$
and that every variable $x_i$ appears in $s$ exactly once. For example, $(x_1+(x_2+x_3))$ and $((x_3+x_1)+x_2)$ are  additive $3$-monomials. We show that additive $n$-monomials are mutually f-equivalent.

\begin{prop}\label{prop_simAddMono}
Let $n\in\N$. Any two additive $n$-monomials are f-equivalent.   
\end{prop}
\proof
An additive $n$-monomial $s$ is standard if  
$$
s=(x_1+(x_2+(\ldots+(x_{n-1}+x_n)\ds)))\;.
$$
We use the inductive argument from the proof of Proposition~\ref{prop_StandMonom}, in which we replace $\cdot$ with $+$, and ET6 
and ET8 with ET5 and ET7.
\eproof

We define formal finite sums of terms.

\begin{defi}\label{def_finSum}
Let $n\in\N_0$ and $t_1$, $\ds$, $t_n$ be $n$ terms. Their sum, denoted by $t_1+\ds+t_n$ or by
$$
\sum_{i=1}^n t_i\;,
$$
is for $n\ge1$ the set of all terms obtained by replacing in all $\frac{1}{n}\binom{2n-2}{n-1}\cdot n!$ additive $n$-
monomials for $i\in[n]$ the occurrence of $x_i$ with $t_i$. For $n=0\in\N_0$ we 
define the sum as $\{0\}\sus T$.
\end{defi}
For example,
$$
\sum_{i=1}^2 (x_i\cdot x_i)=\{((x_1\cdot x_1)+(x_2\cdot x_2)),\,((x_2\cdot x_2)+(x_1\cdot x_1))\}\;.
$$
We determine values of the map $\Psi$ (Definition~\ref{def_psi}) on sums.

\begin{prop}\label{prop_PsiOfSums}
Let $R$ be a~ring, $n\in\N_0$ and $t_1$, $\ds$, $t_n$ be $n$ terms. If $s\in\sum_{i=1}^n t_i$ then
$$
\Psi[s]=\left\{
\begin{array}{lll}
     \emptyset&\ds&n=0\,\text{ and}  \\
     \big(\bigoplus_s\big)_{i=1}^n\Psi[t_i]&\ds&n\ge1\;, 
\end{array}
\right.
$$
where $\oplus_s$ is the addition in $R[x_1,x_2,\ds]$.
\end{prop}
\proof
Let $s$ be as stated. For $n=0$ we have that $s=0$ and $\Psi[s]=\emptyset$ by Definition~\ref{def_psi}. For 
$n=1$ we have that $s=t_1$ and $\Psi[s]=\Psi[t_1]$ holds trivially. For $n>1$ we proceed by induction on $n$. Let 
$$
s_1=(x_1+(x_2+(\ldots+(x_{n-1}+x_n)\ds)))
$$
be the standard additive $n$-monomial and $s'$ be the term obtained from $s_1$ by 
replacing for $i\in[n]$ the occurrence of $x_i$ with $t_i$. Then $\Psi[s]=\Psi[s']$ 
by Propositions~\ref{prop_feqImplseqPol} and \ref{prop_subsXi}. We set $u_i:=t_{i+1}$ and denote by $s''$ the term obtained by replacing in
$$
(x_1+(x_2+(\ldots+(x_{n-2}+x_{n-1})\ds)))
$$
$x_i$ with $u_i$. Using induction for $s''$ we get that $\Psi[s]$ equals to
$${\textstyle
\Psi[s']=\Psi[t_1]\oplus_s\Psi[s'']=
\Psi[t_1]\oplus_s\big(\bigoplus_s\big)_{i=2}^n\Psi[t_i]=\big(\bigoplus_s\big)_{i=1}^n\Psi[t_i]\;.
}
$$
\eproof

We introduce standard terms and their carriers. 

\begin{defi}\label{def_standTerm}
Let $R$ be a~ring. A~term $t$ is standard if
$$
t\in\sum_{i=1}^N (c_{r_i}\cdot\mu_i)
$$
where $N\in\N_0$, $r_i\in R\setminus\{0_R\}$ and the $\mu_i$ are 
$\overline{m(i)}$-monomials with distinct types (so all $N$ tuples 
$\overline{m(i)}\in\N_0^{\infty}$ are mutually distinct). The carrier $S[t]$ of $t$ is the standard polynomial
$$
S[t]:=\{(\overline{m(i)},\,r_i)\;|\;i\in[N]\}\in R[x_1,\,x_2,\,\ds]\;.
$$
\end{defi}
The carrier of any standard term is determined uniquely. For $N=0$ we have $t=0$ and $S[t]=\emptyset$. For $\overline{m(i)}=\emptyset$ we have $\mu_i=1$.

We show that every term is f-equivalent to a~standard term.

\begin{thm}\label{thm_redToStand}
For every term $t$ there is a~standard term $u$ such that $u\stackrel{f}{\approx}t$.    
\end{thm}
\proof
Let $t\in T_{a_R}$ for a~ring $R$. Terms of the form $(v+v')$ are called additive; multiplicative terms are defined 
similarly. A~subterm $v$ of $t$ is  called a~d-subterm if $v=(v_1\cdot v_2)$ with at 
least one $v_i$ additive. A~d-subterm $v$ is deepest if it has the 
maximum depth $\mathrm{dep}(v)$ (Definition~\ref{def_depth}) among all 
d-subterms of $t$. If $t$ has at 
least one d-subterm, in the reduction pair 
$$
\mathrm{red}(t):=(\mathrm{dep}(v),\,k)\in\N^2
$$ 
$v$ is any of the deepest d-subterms and $k$ is their number. If $t$ has no 
d-subterm, we set $\mathrm{red}(t):=0$.

Let $\mathrm{red}(t)=(j,k)$ and $v$ be a~deepest d-subterm of $t$. We may assume 
(by applying ET6) that $v=(v_1\cdot (v_3+v_4))$. We get the term $t'\stackrel{9}{\sim}t$ by replacing in $t$ the 
subterm $v$ with the subterm $v'=((v_1\cdot v_3)+(v_1\cdot v_4))$. The only 
subterms of $t$ affected by this are $v$ and its subterms. It follows that 
$$
\text{$\mathrm{red}(t')=0$ or $\mathrm{red}(t')=(j,\,k-1)$ or $\mathrm{red}(t')=(j',\,k')$ with $j'<j$}\;.
$$
Since $(\N_0,<)$ is a~well ordering, repeating this reduction we eliminate 
all d-subterms in $t$ and obtain a~term $u_0$ such that $u_0\stackrel{f}{\approx}t$
and $u_0$ has no d-subterm.

It follows, by induction on the structure of $u_0$, that
$$
u_0\in\sum_{i=1}^N u_i
$$
where $N\in\N$ and every term $u_i$ has only atomic or multiplicative subterms.

We consider the term $u_1$, other terms $u_i$ are dealt with similarly. Using 
Propositions~\ref{prop_subsXi} and \ref{prop_StandMonom} we bring $u_1$ by  to the form
$$
u_1\stackrel{f}{\approx}(v_1\cdot(v_2\cdot(\ldots\cdot(v_{k-1}\cdot v_k)\ds)))
$$
where $k\in\N$, $v_1$, $\ds$, $v_l$ with $0\le l\le k$ are constants and $v_{l+1}$, 
$\ds$, $v_k$ are variables. Using elementary transformations, we further bring $u_1$ to the form
$$
u_1\stackrel{f}{\approx}(c_r\cdot\mu)
$$
where $r\in R$ and $\mu$ is a~monomial. 

Thus we bring $u_0$ to the form
$$
u_0\in\sum_{i=1}^N(c_{r_i}\cdot\mu_i)
$$
where each $r_i\in R$ and each $\mu_i$ is a~monomial. Using 
Propositions~\ref{prop_subsEquiv}, \ref{prop_StandMonom} and \ref{prop_simAddMono}, and transformations
ET6 and ET9 we bring $u_0$ to the same form
$$
u_0\in\sum_{i=1}^{N'}(c_{r_i'}\cdot\mu_i')
$$
where additionally the monomials $\mu_i'$ have distinct types.

If every $r_i'\ne0_R$, we are done as $u_0\stackrel{f}{\approx}t$ and $u_0$ is 
standard. Else we still have to get rid of the zero coefficients. But this is 
easy to do using elementary transformations. So we finally bring $u_0$ by 
f-equivalence to the form
$$
u\in\sum_{i=1}^{N''}(c_{r_i''}\cdot\mu_i'')
$$
where $N''\in\N_0$,  each $r_i''\in R\setminus\{0_R\}$ and the $\mu_i''$ are monomials with distinct types. This is a~standard term f-equivalent to $t$ and we are done.
\eproof

In the last three propositions we show that the implication $\Rightarrow$ which
we are proving holds for standard terms.

\begin{prop}\label{prop_redToStand1}
Let $R$ be a~ring. Any two standard terms $t$ and $u$ with the same carrier $S[t]=S[u]=:S$ are f-equivalent.
\end{prop}
\proof
Let
$$
S=\{(\overline{m(i)},\,r_i)\;|\;i\in[N]\}
$$
where $N\in\N_0$, $r_i\in R\setminus\{0_R\}$ and the tuples $\overline{m(i)}\in\N_0^{\infty}$ are mutually distinct. Then
$$
t\in\sum_{i=1}^N (c_{r_i}\cdot\mu_i)\,\text{ and }\,u\in\sum_{i=1}^N(c_{r_i}\cdot\nu_i)
$$
where for $i\in[N]$ the terms $\mu_i$ and $\nu_i$ are $\overline{m(i)}$-monomials. By 
Proposition~\ref{prop_StandMonom}, $\mu_i\stackrel{f}{\approx}\nu_i$ for every $i\in[N]$. By 
Proposition~\ref{prop_subsEquiv}, $(c_{r_i}\cdot\mu_i)\stackrel{f}{\approx}
(c_{r_i}\cdot\nu_i)$ for every $i\in[N]$. Thus $t\stackrel{f}{\approx}u$ by 
Propositions~\ref{prop_simAddMono} and \ref{prop_subsXi}.
\eproof

\begin{prop}\label{prop_carrPhi}
Let $R$ be a~ring. For every standard term $t$, 
$$
\Psi[t]=S[t]\in R[x_1,\,x_2,\,\ds]
$$
where $S[t]$ is the carrier of $t$ (Definition~\ref{def_standTerm}) and $\Psi$ is the map of Definition~\ref{def_psi}.
\end{prop}
\proof
We proceed by induction on the complexity of
$$
S[t]=\{(\overline{m(i)},\,r_i)\;|\;i\in[N]\},\ N\in\N_0\;.
$$
If $N=0$ then $S[t]=\emptyset$, $t=0$ and $\Psi[t]=\emptyset=S[t]$. If $N\ge1$ then 
$$
t\in\sum_{i=1}^N(c_{r_i}\cdot\mu_i)
$$
where each $\mu_i$ is an $\overline{m(i)}$-term. Then by Propositions~\ref{prop_PsiOfSums} and \ref{prop_onMonom2},
$${\textstyle
\Psi[t]=\big(\bigoplus_s\big)_{i=1}^N\Psi[(c_{r_i}\cdot\mu_i)]\stackrel{s}{=}
\big(\bigoplus_s\big)_{i=1}^N r_i x^{\overline{m(i)}}=S[t]
}
$$
(with addition $\oplus_s$ in $R[x_1,x_2,\ds]$).
\eproof

\begin{prop}\label{prop_finalStep}
If $t$ and $u$ are s-equivalent standard terms, then $t$ and $u$ are f-equivalent.    
\end{prop}
\proof
This follows from Definition~\ref{def_sEquiv} and Propositions~\ref{prop_redToStand1} and \ref{prop_carrPhi}.
\eproof

We finish the proof of the implication $\Rightarrow$. Let $t$ and $u$ be s-equivalent terms. By 
Theorem~\ref{thm_redToStand} there are standard terms $t'$ and $u'$ such that 
$t\stackrel{f}{\approx}t'$ and $u\stackrel{f}{\approx}u'$. By 
Proposition~\ref{prop_feqImplseqPol}, also $t\stackrel{s}{\approx}t'$ and 
$u\stackrel{s}{\approx}u'$. The transitivity of $\stackrel{s}{\approx}$ gives that $t'\stackrel{s}{\approx}u'$. By 
Proposition~\ref{prop_finalStep},  $t'\stackrel{f}{\approx}u'$. The transitivity of $\stackrel{f}{\approx}$ gives that $t\stackrel{f}{\approx}u$. {\bf The proof of 
Theorem~\ref{thm_sEqIsfEq} is complete}.
\eproof

\noindent
This concludes the initial part of our article on formal polynomials. In the rest we turn to Boolean terms and to independence of events.

\section[Boolean algebras and Boolean formulas]{Boolean algebras and Boolean formulas}\label{sec_BooTermPoly}

In this section we review Boolean algebras and Boolean terms (formulas), and 
introduce evaluation maps. 

In the definition of Boolean algebras we follow \cite[Chapter 11.3]{macl_birk}.
A~{\em Boolean algebra} 
$$
A=(A,\,\vee,\,\wedge,\,\neg,\,0_A,\,1_A)
$$
is a~sixtuple of a~base set $A\ne\emptyset$, binary operations
$$
\vee,\,\wedge\cc A\times A\to A
$$
on it called respectively the {\em join} and the {\em meet}, a~unary operation
$$
\neg\cc A\to A
$$
on it called the {\em complement} and distinct constants $0_A,1_A\in A$ 
called respectively {\em zero} and {\em one}. To agree with Boolean formulas, 
we denote the complement by $\neg$ and return to the more common notation ${}^c$ later. These operations satisfy 
the following Boolean axioms. For every $a\in A$ it holds that $0_A\vee a=a$, 
$0_A\wedge a=0_A$, $1_A\vee a=1_A$, $1_A\wedge a=a$,
$$
a\vee \neg a=1_A\;\text{ and }\;a\wedge \neg a=0_A\;.
$$
Both operations $\vee$ and $\wedge$ are commutative, associative, idempotent ($a\vee a=a\wedge a=a$) and mutually distributive,
$$
a\wedge(b\vee c)=(a\wedge b)\vee(a\wedge c)\,\text{ and }\,a\vee(b\wedge c)=(a\vee b)\wedge(a\vee c)\;.
$$
Finally, for every $a,b\in A$ the absorption rules 
$$
a\vee(a\wedge b)=a\wedge(a\vee b)=a
$$ 
hold.

By the commutativity and associativity of $\vee$ and $\wedge$, for any nonempty finite set $B\sus A$ the results of the
$|B|$-fold join and meet
$$
{\textstyle
\bigvee_{b\in B}b\in A\,\text{ and }\,
\bigwedge_{b\in B}b\in A
}
$$
are uniquely determined. In more detail, if $B=\{b_1,\ds,b_n\}\ne\emptyset$ with $|B|=n$ and $\circ\in\{\vee,\wedge\}$, then for any permutation $\pi$ of $[n]$ and any proper binary
bracketing of
$$
b_{\pi(1)}\,\circ\,b_{\pi(2)}\,\circ\,b_{\pi(3)}\,\circ\,\ds\,\circ\,b_{\pi(n)}\;,
$$
the resulting element of $A$ is independent of $\pi$ and the bracketing. 
In a~proper binary bracketing we place in any of the $\frac{1}{n}\binom{2n-2}{n-1}$ 
possible ways $n-1$ pairs of brackets $($ and $)$ before, between and after the 
above displayed $2n-1$ symbols so that the bracketing is proper and every pair of 
associated brackets encloses two arguments, $(\ds\,\circ
\,\ds)$. It follows more formally (i.e., more precisely) from Propositions~\ref{prop_feqImplseq}, \ref{prop_longDisj} and \ref{prop_longConj} below. For example,
for $n=3$ we have twelve possibilities
$$
((b_{\pi(1)}\,\circ\,b_{\pi(2)})\,\circ\,b_{\pi(3)})\,\text{ and }\,
(b_{\pi(1)}\,\circ\,(b_{\pi(2)}\,\circ\,b_{\pi(3)}))\;.
$$

We call elements $a_1$, $\ds$, $a_n$ in a~Boolean algebra $A$ {\em (mutually) 
disjoint} if $a_i\wedge a_j=0_A$ whenever $i\ne j$. 

We consider the functional symbols $F:=\{\vee,\wedge,\neg,0,1\}$ and the arity 
function $a\cc F\to\N_0$,
$$
a[\vee]=a[\wedge]=2,\ a[\neg]=1\,\text{ and }\,a[0]=a[1]=0\;.
$$
We use the same symbols $\vee$, $\wedge$ and $\neg$
for operations in different Boolean algebras and for formal functional symbols. Any Boolean algebra has arity $a$. We do not develop Boolean version of Theorem~\ref{thm_sEqIsfEq} and instead concentrate on the Boolean Independence Theorem. Thus there are no $A$-constants in $F$.

We fix notation and terminology for $a$-terms.

\begin{defi}\label{def_Bterms}
Let $n\in\N$ and $a$ be the above arity. In the rest of the article we call elements in the sets
$$
T:=T_a\,\text{ and }\,T[x_1,\,\ds,\,x_n]:=T_a[x_1,\,\ds,\,x_n]
$$
terms or Boolean formulas.
\end{defi}
For instance,
$$
t=\neg(((0\vee x_4)\wedge \neg(x_5)))
$$
is a~term (or a~Boolean formula). It lies in $T$ and in any $T[x_1,\ds,x_n]$ for $n\ge 5$. Variables and their complements, $x_i$ and $\neg(x_i)$ for $i\in\N$, are called {\em literals}.

We introduce evaluation maps.

\begin{defi}\label{def_evalMBool}
Let $n\in\N$, let $A=(A,\vee,\wedge,\neg,0_A,1_A)$ be a~Boolean algebra and $t$ in $T[x_1,\ds,x_n]$ be a~term. The evaluation map
$$
\Phi_A[t;\,y_1,\,\ds,\,y_n]=\Phi_{\mathcal{A}}[t;\,y_1,\,\ds,\,y_n]\cc A^n\to A
$$
associated to $t$ is a~particular case of the function $\Phi_{\mathcal{A}}$ in Definition~\ref{def_evalMap}, with the algebra $\mathcal{A}$ of arity $a$ being $A$.
\end{defi}
If $n$ is minimum such that $t\in T_a[x_1,\ds,x_n]$ and $m>n$, then the values of $\Phi_A[t;\,y_1,\ds,y_m]$ are independent of the arguments $y_i$ with $i>n$.

\section[DNF terms]{DNF terms}\label{sec_DNFterms}
We parallel Section~\ref{sec_polynomials} to some extent. In Definition~\ref{def_BoolTrans} to 
Proposition~\ref{prop_subsEquiv2} we introduce Boolean transformations of 
terms, f-equivalence of terms and evaluation maps. In 
Propositions~\ref{prop_longDisj}--\ref{prop_01areDNF} we introduce DNF terms. In Theorem~\ref{thm_DNF} we prove 
that every term is f-equivalent to a~DNF term. 

\begin{defi}\label{def_BoolTrans}
Two terms $t$ and $t'$ are 
(Boolean) transforms of one another, written $t\sim t'$ and $t'\sim t$, if for some 
$u,u',u''\in T$ the term $t'$ arises from the term $t$ by replacing a~subterm $v$ of $t$ with the subterm $v'$
by one of the nineteen rules below.
\begin{description}
    \item[{\rm BT1}] $v=(0\wedge u)$ and $v'=0$.
    \item[{\rm BT2}] $v=(0\vee u)$ and $v'=u$.
    \item[{\rm BT3}] $v=(1\wedge u)$ and $v'=u$.
    \item[{\rm BT4}] $v=(1\vee u)$ and $v'=1$.
    \item[{\rm BT5}] $v=(u\wedge \neg(u))$ and $v'=0$.
    \item[{\rm BT6}] $v=(u\vee \neg(u))$ and $v'=1$.
    \item[{\rm BT7}] $v=(u\vee u')$ and $v'=(u'\vee u)$.
    \item[{\rm BT8}] $v=(u\wedge u')$ and $v'=(u'\wedge u)$.
    \item[{\rm BT9}] $v=(u\vee (u'\vee u''))$ and $v'=((u\vee u')\vee u'')$.
    \item[{\rm BT10}] $v=(u\wedge (u'\wedge u''))$ and $v'=((u\wedge u')\wedge u'')$.
    \item[{\rm BT11}] $v=(u\vee u)$ and $v'=u$.
    \item[{\rm BT12}] $v=(u\wedge u)$ and $v'=u$.
    \item[{\rm BT13}] $v=(u\vee (u'\wedge u'')$ and $v'=((u\vee u')\wedge (u\vee u''))$.
    \item[{\rm BT14}] $v=(u\wedge (u'\vee u''))$ and $v'=((u\wedge u')\vee(u\wedge u''))$.
    \item[{\rm BT15}] $v=u\vee(u\wedge u')$ and $v'=u$.
    \item[{\rm BT16}] $v=u\wedge(u\vee u')$ and $v'=u$.
    \item[{\rm BT17}] $v=\neg(\neg(u))$ and $v'=u$.
    \item[{\rm BT18}] $v=\neg((u\vee u'))$ and $v'=(\neg(u)\wedge \neg(u'))$.
    \item[{\rm BT19}] $v=\neg((u\wedge u'))$ and $v'=(\neg(u)\vee \neg(u'))$.
\end{description}
\end{defi}
The relation $\sim$ is obviously symmetric. Notation $t\stackrel{i}{\sim}t'$ means that 
$t\sim t'$ by $\mathrm{BT}i$, $i\in[19]$. Rules BT1--BT16 correspond to Boolean 
axioms. Rule BT17 corresponds to the double negation (complement) rule and rules BT18 and 
BT19 correspond to de Morgan rules; their counterparts hold in every Boolean algebra 
(\cite[Chapter~11.4]{macl_birk}).

We introduce f-equivalence of terms.

\begin{defi}\label{def_fekviv2}
Two terms $t$ and $u$ are f-equivalent, written $t\stackrel{f}{\approx}u$, 
if there exist $k+1$ terms $t_0$, $t_1$, $\ds$, $t_k$, $k\in\N$, such that 
$$
t=t_0\sim t_1\sim\ds\sim t_k=u\ (\text{i.e., $t_{i-1}\sim t_i$ for $i\in[k]$})\;.
$$
Thus one can transform $t$ in $u$ by a~series of transformations $\sim$ listed in Definition~\ref{def_BoolTrans}.
\end{defi}
The relation $\stackrel{f}{\approx}$ (cf. Definition~\ref{def_fekviv}) is an equivalence on 
Boolean formulas. For example, $((\neg(x_1)\wedge x_2)\wedge x_1)\stackrel{f}{\approx} 0$ because
$$
((\neg(x_1)\wedge x_2)\wedge x_1)\stackrel{8,\,10}{\sim}((x_1\wedge \neg(x_1))\wedge x_2)
\stackrel{5}{\sim}(0\wedge x_2)\stackrel{1}{\sim}0\;.
$$
Also, $\neg(0)\stackrel{f}{\approx}1$ and $\neg(1)\stackrel{f}{\approx}0$ because
$$
\neg(0)\stackrel{2}{\sim}(0\vee\neg(0))\stackrel{6}{\sim}1\,\text{ and }\,
\neg(1)\stackrel{3}{\sim}(1\wedge\neg(1))\stackrel{5}{\sim}0\;.
$$

Recall the map $\Phi_A$ of Definition~\ref{def_evalMBool}.
The next proposition is a~Boolean analogue of  
Proposition~\ref{prop_feqImplseqPol}. 

\begin{prop}\label{prop_feqImplseq}
For every Boolean algebra $A$ and every two f-equivalent terms $t$ and $u$ in $T[x_1,\ds,x_n]$ one has that
$$
\Phi_A[t,\;y_1,\,\ds,\,y_n]=\Phi_A[u,\;y_1,\,\ds,\,y_n]\;.
$$
\end{prop}
\proof
The proof is similar to that of Proposition~\ref{prop_feqImplseqPol}.
\eproof

We proceed to the definition of DNF terms and first again note that substitution in subterms preserves 
f-equivalence.

\begin{prop}\label{prop_subsEquiv2}
If a~term $t'$ arises from another term $t$ by replacing a~subterm $u$ of $t$ with an f-equivalent subterm $u'$ then 
$t\stackrel{f}{\approx}t'$.
\end{prop}
\proof
The proof  is similar to that of Proposition~\ref{prop_subsEquiv}
\eproof

Let $t_1$, $\ds$, $t_n$ be terms, $n\in\N_0$. Their {\em disjunction} is the set of terms denoted by 
$\bigvee_{i=1}^n t_i$ and obtained as follows.  For $n\ge1$ we take all terms 
$$
u\in\big(\{\,),\ (,\ \vee\,\}\cup\{x_1,\,\ds,\,x_n\}\big)^*
$$
in which every variable $x_i$ appears exactly once and in every $u$ 
we replace for $i\in[n]$ the occurrence of $x_i$ with the term $t_i$. 
For $n=0$ we set $\bigvee_{i=1}^0 t_i:=\{0\}$. For example,
$$
\text{$\bigvee_{i=1}^3 t_i$ is the $12$-element set $\{((t_1\vee t_2)\vee t_3),\,\ds,\,(t_3\vee(t_2\vee t_1))\}$}\;.
$$

We show that the terms in the disjunction of $t_1$, $\ds$, $t_n$ are mutually f-equivalent.

\begin{prop}\label{prop_longDisj}
For any terms $t_1$, $\ds$, $t_n$, $n\in\N_0$, any two terms in $\bigvee_{i=1}^nt_i$ are 
f-equivalent.
\end{prop}
\proof
We proceed as in the proof of Proposition~\ref{prop_StandMonom} and show that every term in the disjunction of $t_1$, $\ds$, $t_n$ is f-equivalent to the standard form
$$
(t_1\vee(t_2\vee(\ds\vee(t_{n-1}\vee t_n)\ds)))\;.
$$
\eproof

The {\em conjunction} $\bigwedge_{i=1}^n t_i$ of the terms $t_1$, $\ds$, $t_n$, $n\in\N_0$, is defined similarly to the
disjunction, with $\vee$ replaced by $\wedge$. For $n=0$ we set $\bigwedge_{i=1}^0 t_i:=\{1\}$ This replacement turns the previous proof in the proof of the next proposition.

\begin{prop}\label{prop_longConj}
For any terms $t_1$, $\ds$, $t_n$ with $n\in\N_0$ the conjunctions in $\bigwedge_{i=1}^n t_i$ are mutually 
f-equivalent.    
\end{prop}

Let $n\in\N_0$, 
$X=\{m_1<m_2<\ds<m_n\}\sus\N$ and for $i\in[n]$, $t_i:=x_{m_i}$ or 
$t_i:=\neg(x_{m_i})$. If $t\in\bigwedge_{i=1}^n t_i$, we
say that $t$ is an $X$-$\mathrm{monomial}_{\wedge}$. The vector $\tau(t)=(c_1,\ds,c_n)\in\{0,1\}^n$,
where $c_i=0$ if $t_i=x_{m_i}$ and  $c_i=1$ if $t_i=\neg(x_{m_i})$, is the 
{\em type of $t$}. For $n=0$ we set $\tau(t)=\emptyset$. We can simplify 
conjunctions and disjunctions of monomials as follows.

\begin{prop}\label{prop_diffTypes0}
If $\mu$ and $\nu$ are $X$-$\mathrm{monomials}_{\wedge}$ with the same type then
$$
(\mu\vee\nu)\stackrel{f}{\approx}\mu\stackrel{f}{\approx}\nu\;.
$$
If $\mu$ and $\nu$ are $X$-$\mathrm{monomials}_{\wedge}$ with different types then
$$
(\mu\wedge\nu)\stackrel{f}{\approx}0\;.
$$
\end{prop}
\proof
Suppose that $\mu$ and $\nu$ are $X$-$\mathrm{monomials}_{\wedge}$ with the 
same type. By Proposition~\ref{prop_longConj}, $\mu\stackrel{f}{\approx}\nu$. Thus by 
Proposition~\ref{prop_subsEquiv2} also $\mu\stackrel{11}{\sim}(\mu\vee\mu)\stackrel{f}{\approx}(\mu\vee\nu)$.

Suppose that $\mu$ and $\nu$ have different types. By Proposition~\ref{prop_longConj},
$$
\mu\stackrel{f}{\approx}(x_i\wedge\mu')\,\text{ and }\nu\stackrel{f}{\approx}(\neg(x_i)\wedge\nu')
$$
for some $i\in X$ and some $(X\setminus\{i\})$-$\mathrm{monomials}_{\wedge}$ $\mu'$ and $\nu'$, where $x_i$ and $\neg(x_i)$ may be swapped. Thus
\begin{eqnarray*}
&&(\mu\wedge\nu)\stackrel{f}{\approx}((x_i\wedge\mu')\wedge(\neg(x_i)\wedge\nu'))\stackrel{f}{\approx}\\
&&\stackrel{f}{\approx}
((x_i\wedge\neg(x_i))\wedge(\mu'\wedge\nu'))
\stackrel{f}{\approx}0\;,
\end{eqnarray*}
by Proposition~\ref{prop_subsEquiv2}, BT1, BT5, BT8 and BT10.
\eproof 

\begin{defi}\label{def_DNFandCNF}
Let $n\in\N$. A~term $t$ is an $n$-standard {\em DNF} (disjunctive normal form) term if
$$
t\in\bigvee_{i=1}^k\mu_i
$$
where $k\in\N_0$ and the $\mu_i$ are $[n]$-$\mathrm{monomials}_{\wedge}$ with mutually distinct types.
\end{defi}
For example,
$$
(((x_2\wedge x_3)\wedge\neg(x_1))\vee
(x_1\wedge(x_3\wedge\neg(x_2))))
$$
is a~3-standard DNF term.

\begin{prop}\label{prop_01areDNF}
For every $n\in\N$ the constants $0$ and $1$ are f-equivalent to $n$-standard {\em DNF} terms because $0\in\bigvee_{i=1}^0\mu_i$ and 
$1\stackrel{f}{\approx}t$ for every
$$
t\in\bigvee_{i=1}^{{2^n}}\mu_i
$$
and every selection of $2^n$ $[n]$-$\mathrm{monomials}_{\wedge}$ $\mu_1,\ds,\mu_{2^n}$ with mutually distinct types (so each type appears exactly once). 
\end{prop}
\proof
The first claim is just the definition
of the empty disjunction. We prove the second f-equivalence. By 
Propositions~\ref{prop_longDisj} and \ref{prop_longConj} it suffices to show 
that for one term $t$ in the displayed disjunction for some selection of $\mu_1$, 
$\ds$, $\mu_{2^n}$ one has that $t\stackrel{f}{\approx}1$. We set $t_i:=
(x_i\vee\neg(x_i))\stackrel{6}{\sim}1$ and take 
$$
t:=(t_1\wedge(t_2\wedge(\ds\wedge(t_{n-1}\wedge t_n))))\;.
$$
On the one hand we have $t\stackrel{f}{\approx}1$ by 
Proposition~\ref{prop_subsEquiv2} and BT3. On the other hand it follows by repeated 
application of BT14 and BT8 that $t$ is f-equivalent to a~term in $\bigvee_{i=1}^{{2^n}}\mu_i$ for some $[n]$-$\mathrm{monomials}_{\wedge}$ $\mu_i$
with distinct types.
\eproof

We show that every term is f-equivalent to a~DNF term. The reader may wish to compare the proof 
below with the proof of Theorem~\ref{thm_redToStand}.

\begin{thm}\label{thm_DNF}
For every term $t\in T[x_1,\ds,x_n]$, $n\in\N$, there is an $n$-standard {\em 
DNF} term $u$ such that $u\stackrel{f}{\approx}t$.
\end{thm}
\proof
Suppose that $t$ is a~term in $T[x_1,\ds,x_n]$ with $n\in\N$. A~subterm $v$ of $t$ is called an n-subterm if  
$v=\neg(v_1)$ and $v_1$ is not a~variable. This subterm is 
called deepest if it has the maximum depth $\mathrm{dep}(v)$ among all n-subterms of $t$. If $t$ has
at least one n-subterm, in the reduction pair
$$
\mathrm{red}(t):=(\mathrm{dep}(v),\,k)\in\N^2
$$
$v$ is any of the deepest n-subterms of $t$ and $k$ is their number. Else we set $\mathrm{red}(t):=0$. 

Let $\mathrm{red}(t)=(j,k)$ and $v=\neg(v_1)$ be a~deepest n-subterm of $t$. Then (i) $v_1=\neg(v_2)$ or 
(ii) $v_1=(v_2\vee v_3)$ or (iii) $v_1=(v_2\wedge v_3)$ or (iv) $v_1=0$ or (v) $v_1=1$. We get the 
term $t'\stackrel{f}{\approx}t$ by replacing $v$ with the subterm $v'=v_2$ via BT17 in case (i), 
with $v'=(\neg(v_2)\wedge\neg(v_3))$ via BT18 in case (ii), with $v'=(\neg(v_2)\vee\neg(v_3))$ via BT19 in case 
(iii) and with $v'=1$, resp. $v'=0$, via the above f-equivalences for the complement 
of $0$, resp. $1$, in case (iv), resp. case (v). As before it follows that
$$
\text{$\mathrm{red}(t')=0$ or $\mathrm{red}(t')=(j,k-1)$ or $\mathrm{red}(t')=(j',k')$ with $j'<j$}\;.
$$

Repeating this reduction we get a~term 
$u_0\stackrel{f}{\approx}t$ with 
$\mathrm{red}(u_0)=0$, i.e. every subterm of $u_0$ is a~disjunction $(v\vee v')$, a~conjunction $(v\wedge v')$, a~literal or 
a~constant ($0$ or $1$). Using BT1--BT4 we can eliminate the last possibility, except when $u_0=0$ 
or $u_0=1$. Thus we get a~term $u_0\stackrel{f}{\approx}t$ such that 
$$
\text{$u_0=0$ 
or $u_0=1$ or $u_0\in\big(\{(,\,),\,\vee,\,\wedge\}\cup\{x_1,\,\neg(x_1),\,\ds,\,x_n,\,\neg(x_n)\}\big)^*=:A^*$}\;.
$$
In the first two cases we are done since $0$ and $1$ are f-equivalent to $n$-standard DNF terms by 
Proposition~\ref{prop_01areDNF}. We deal with the third case.

Let $u_0$ be a~term in $A^*$. A~subterm $v$ of $u_0$ is called a~d-subterm if (i) $v=(v_1\wedge(v_2\vee v_3))$ or (ii) $v=((v_2\vee v_3)\wedge v_1)$. We define deepest d-subterms of $u_0$ in the obvious way and if $u_0$
has a~d-subterm, we set as before
$$
\mathrm{red}(u_0):=(\mathrm{dep}(v),\,k)\in\N^2
$$
where $v$ is one of the deepest d-subterms of $u_0$
and $k$ is their number. Else we set $\mathrm{red}(u_0):=0$. Let $\mathrm{red}(u_0)=(j,k)$ and $v$
be a~deepest d-subterms of $u_0$. We obtain the term $u_1\stackrel{f}{\approx}u_0$ by replacing $v$
with the subterm $v'=((v_1\wedge v_2)\vee(v_1\wedge v_3))$ via BT14 in case (i) and via BT8 and BT14 in 
case (ii). Again it follows that
$$
\text{$\mathrm{red}(u_1)=0$ or $\mathrm{red}(u_1)=(j,k-1)$ or $\mathrm{red}(u_1)=(j',k')$ with $j'<j$}\;.
$$

Repeating this reduction we get a~term $u_1\stackrel{f}{\approx}t$ in 
$A^*$ with no d-subterm. We see by induction on the structure of $u_1$ 
that
$$
u_1\in\bigvee_{i=1}^k\mu_i
$$
where $k\in\N$ and each $\mu_i$ is a~(nonempty) conjunction of some, possibly repeating, literals $x_1$, $\ds$, $\neg(x_n)$. 

In the final phase of the proof we transform the $\mu_i$ in $[n]$-$\mathrm{monomials}_{\wedge}$ with 
distinct types. We consider $\mu_1$, the other conjunctions $\mu_i$ are treated similarly. Using 
Proposition~\ref{prop_longConj} and BT12 we remove in $\mu_1$ all duplicate literals. If for some 
$i\in[n]$ both $x_i$ and $\neg(x_i)$ are disjoint subterms of (the reduced) $\mu_1$, then by using 
Proposition~\ref{prop_longConj}, BT5 and BT1 we get that $\mu_1\stackrel{f}{\approx}0$. By
Proposition~\ref{prop_longDisj} and BT2 we may omit these zero disjuncts in $u_1$. Thus we get a~term
$u_2\stackrel{f}{\approx}t$ such that
$$
u_2\in\bigvee_{i=1}^{k'}\mu_i'
$$
where $k'\in\N_0$ and the $\mu_i'$ are $X_i$-$\mathrm{monomials}_{\wedge}$ for nonempty sets $X_i\sus[n]$. If $k'=0\in\N_0$ then $u_2=0\in T$ and, as we know, we are done.

Suppose that $k'>0$. It remains to extend each $X_i$ to $[n]$ and then to remove monomials with 
duplicate types. If there is an index $j\in[n]\setminus X_1$ then
$$
\mu_1'\stackrel{3}{\sim}1\wedge\mu_1'
\stackrel{6}{\sim}(x_j\vee\neg(x_j))\wedge\mu_1'
\stackrel{8,\,9}{\sim}((\mu_1'\wedge x_j)
\vee(\mu_1'\wedge\neg(x_j)))
$$
and using Proposition~\ref{prop_subsEquiv2} we can replace $\mu_1'$ with the term $(\nu_1\vee\nu_1')$
where $\nu_1$ and $\nu_1'$ are $X_1\cup\{j\}$-$\mathrm{monomials}_{\wedge}$. Repeating 
these extensions we get a~term $u_3\stackrel{f}{\approx}t$ such that
$$
u_3\in\bigvee_{i=1}^{k''}\mu_i''
$$
where $k''\in\N$ and the $\mu_i''$ are $[n]$-$\mathrm{monomials}_{\wedge}$. Finally, the 
monomials $\mu_i''$ with duplicate types may be omitted in $u_3$ due to Propositions~\ref{prop_subsEquiv2}, \ref{prop_longDisj} and \ref{prop_diffTypes0}.
\eproof

\noindent
The previous proof is a~more formal and more precise version of the proof in \cite[Chapter 
11.5]{macl_birk}.

If $S$ and $S'$ are nonempty finite sets of terms, we define the nonempty set
$$
S\wedge S':=\{(t\wedge u)\;|\;t\in S,\,u\in S'\}\sus T\;.
$$
In the proof of BIT we need the following result.

\begin{prop}\label{prop_longDistr}
For every terms $t_1$, $\ds$, $t_m$, $u_1$, $\ds$, $u_n$ with $m,n\in\N_0$
it is true that
$$
\bigvee_{i=1}^m t_i\wedge\bigvee_{j=1}^n u_i\stackrel{f}{\approx}
\bigvee_{i,\,j=1}^{m,\,n}(t_i\wedge u_j)
$$
--- every term on the left side of  $\stackrel{f}{\approx}$ is f-equivalent to every term on the right side.
\end{prop}
\proof
For $m=0$ or $n=0$ it holds, both sides are $\{0\}$. Let $m,n>0$.
We show that all terms on both sides of $\stackrel{f}{\approx}$ are mutually f-equivalent. It follows from 
Propositions~\ref{prop_subsEquiv2} and \ref{prop_longDisj} that any two terms on the 
same side are f-equivalent. It therefore suffices if we show that any term on 
the left side is f-equivalent to a~term on the right side. For $m=1$ it holds by repeated application of BT14. We assume that $m\ge2$ and proceed by induction on $m$. Let 
$t\in\bigvee_{i=1}^m t_i$ and $u\in\bigvee_{j=1}^n u_j$ be arbitrary. Then $t=(t'\vee t'')$ for some $t'\in\bigvee_{i\in X}t_i$ and $t''\in\bigvee_{i\in Y}t_i$ where the sets $X$ and $Y$ partition $[m]$. We have that 
$$
(t\wedge u)=((t'\vee t'')\wedge u)
\stackrel{8,\,14}{\sim}((t'\wedge u)\vee(t''\wedge u))\;.
$$
The initial conjunction is an arbitrary term on the left side. By induction, 
$$
(t'\wedge u)\stackrel{f}{\approx}v\in
\bigvee_{i\in X}\bigvee_{j=1}^n(t_i\wedge u_j)
\,\text{ and }\,
(t''\wedge u)\stackrel{f}{\approx}w\in
\bigvee_{i\in Y}\bigvee_{j=1}^n(t_i\wedge u_j)\;.
$$
By Proposition~\ref{prop_subsEquiv2},
$$
(t\wedge u)\stackrel{f}{\approx}(v\vee w)\in
\bigvee_{i\in X\cup Y}\bigvee_{j=1}^n(t_i\wedge u_j)=\bigvee_{i,\,j=1}^{m,\,n}(t_i\wedge u_j)\;.
$$
\eproof

\section[The Boolean Independence Theorem]{The Boolean Independence Theorem}\label{sec_BoolInde}

We propose a~formalization of the type of probability used in the 
Probabilistic Method in Combinatorics. Problems investigated there usually (but not always) involve only finitely 
many events. The full $\sigma$-additive probability is usually not needed and simpler formalization suffices. 

It is a~bizarre fact that the basic monographs \cite{alon_spen} for the 
Probabilistic Method contain no theoretical foundations;  
there is no definition of independence of events and random variables, of 
probability spaces, of martingales, etc., all of which are basic tools in 
\cite{alon_spen}. Thus it may be useful if we formalize here at least the independence facet. 

\begin{defi}\label{def_probSpace}
A~finitary probability space, abbreviated {\em FPS}, is a~pair
$$
P=(A,\,\mathrm{Pr})
$$
of a~Boolean algebra $A=(A,\vee,\wedge,{}^c,0_A,1_A)$ (from now on we use ${}^c$ to denote the operation of complement) and a~probability map
$$
\mathrm{Pr}\cc A\to[0,\,1]
$$
satisfying that $\mathrm{Pr}(0_A)=0$, $\mathrm{Pr}(1_A)=1$ and for every $a,b\in A$,
$$
\mathrm{Pr}(a\vee b)=\mathrm{Pr}(a)+\mathrm{Pr}(b)-\mathrm{Pr}(a\wedge b)\;.
$$
\end{defi}
Setting $b:=a^c$ we get that $\mathrm{Pr}(a^c)=1-\mathrm{Pr}(a)$ for every $a\in A$.
The elements $a\in A$ are called {\em events} and the values $\mathrm{Pr}(a)$ 
are their {\em probabilities}. 

The FPS used in the Probabilistic Method are often given by a~nonempty finite set 
$M$ with real weights $p_x\in[0,1]$, 
$x\in M$, such that $\sum_{x\in M}p_x=1$. One then has the Boolean algebra $(A,\vee,\wedge,{}^c,0_A,1_A)$ with
$$
A:=\mathcal{P}(M),\,\vee=\cup,\,\wedge=\cap,\,a^c=M\setminus a,\,0_A=\emptyset,\,1_A=M
$$
and the probability ($a\sus M$)
$$
{\textstyle
\mathrm{Pr}(a):=\sum_{x\in a}p_x\;.
}
$$
One often works with an infinite sequence of such FPS.

The following summation formula is often used.

\begin{cor}\label{cor_disjEven}
Suppose that $a_1$, $\ds$, $a_n$, $n\in\N$, are mutually disjoint events in a~{\em FPS} $(A,\mathrm{Pr})$. Then 
$$
{\textstyle
\mathrm{Pr}\big(\bigvee_{j=1}^n a_j\big)=\sum_{j=1}^n\mathrm{Pr}(a_j)\;.
}
$$
\end{cor}
\proof
For $n=1$ the equality holds trivially. For $n\ge2$ we use induction. We set 
$a:=a_1$ and $b:=\bigvee_{i=2}^n a_i$. By repeated application of Boolean axioms and by the definition of the big disjunction $\bigvee$ (more precisely, by Propositions~\ref{prop_feqImplseq} and \ref{prop_longDistr} and by some Boolean transformations) we get that 
$${\textstyle
a\wedge b=\bigvee_{i=2}^n(a_1\wedge a_i)=
\bigvee_{i=2}^n 0_A=0_A\;.
}
$$
Using Definition~\ref{def_probSpace} and 
induction we thus have
\begin{eqnarray*}
{\textstyle
\mathrm{Pr}\big(\bigvee_{j=1}^n a_j\big)}&=&
\mathrm{Pr}(a\vee b)=\mathrm{Pr}(a)+\mathrm{Pr}(b)-\mathrm{Pr}(a\wedge b)\\
&=&\mathrm{Pr}(a_1)+{\textstyle\sum_{j=2}^n\mathrm{Pr}(a_j)-0=
\sum_{j=1}^n\mathrm{Pr}(a_j)}\;.
\end{eqnarray*}
\eproof

We review independence of events. Two events $a$ and $b$ in a~FPS $(A,\mathrm{Pr})$ are {\em 
independent} if $\mathrm{Pr}(a\wedge b)=\mathrm{Pr}(a)\mathrm{Pr}(b)$. More
generally, finitely many events $a_i$, $i\in I$, are {\em independent} if for every subset $X\sus I$ one has that
$$
{\textstyle
\mathrm{Pr}\big(\bigwedge_{i\in X}a_i\big)=\prod_{i\in X}\mathrm{Pr}(a_i)\;.
}
$$
Another generalization is the following.
\begin{defi}\label{def_Indep}
Let $a_i$, $i\in I$, and $b_i$, $i\in J$, be two finite tuples of events in a~{\em FPS} $(A,\mathrm{Pr})$. We say that the former events are independent of the 
latter events if for every subsets $X\sus I$ and $Y\sus J$ one has that
$$
{\textstyle
\mathrm{Pr}\big(\bigwedge_{i\in X}a_i\wedge\bigwedge_{i\in Y}b_i\big)=\mathrm{Pr}\big(\bigwedge_{i\in X}a_i\big)\mathrm{Pr}\big(\bigwedge_{i\in Y}b_i\big)\;.
}
$$
\end{defi}
In particular, an event $a$ is {\em independent of a~finite tuple of events
$b_i$, $i\in J$}, if for every subset $Y\sus J$ one has that
$$
{\textstyle
\mathrm{Pr}\big(a\wedge\bigwedge_{i\in Y}b_i\big)=\mathrm{Pr}(a)\mathrm{Pr}\big(\bigwedge_{i\in Y}b_i\big)\;.
}
$$
It follows that if $a_1$, $\ds$, $a_n$ are independent events and $1\le m\le n$ then the events 
$a_1$, $\ds$, $a_m$ are independent of the events $a_{m+1}$, $\ds$, $a_n$.

Besides Corollary~\ref{cor_disjEven}, another tool in the proof of 
Theorem~\ref{thm_indepThm} is the following proposition. It says that 
complementing events preserves independence.

\begin{prop}\label{prop_indepCompl}
Suppose that $a_1$, $\ds$, $a_{2n}$, $n\in\N$, are events in a~{\em FPS} such 
that the events $a_1$, $\ds$, $a_n$ are independent of the events $a_{n+1}$, 
$\ds$, $a_{2n}$. Then for any selection of events $b_i$, $i\in[2n]$, such that for every $i$ either $b_i=a_i$ or $b_i=a_i^c$, the events $b_1$, $\ds$, $b_n$ are independent of the 
events $b_{n+1}$, $\ds$, $b_{2n}$.
\end{prop}
\proof
Let $a_i$ and $b_i$ be as stated, $X\sus[n]$, $Y\sus\{n+1,\ds,2n\}$ and 
$s\in\N_0$ be the number of events $b_i$
such that $i\in X\cup Y$ and $b_i=a_i^c$.
We proceed by induction on $s$. If $s=0$, the required equality (see below) trivially holds. Let $s>0$ and $j\in X\cup Y$, 
say $j\in X$ (for $j\in Y$ we argue similarly), be such that $b_j=a_j^c$. We consider the events 
$a:=\bigwedge_{i\in X}b_i$, $a':=\bigwedge_{i\in X\setminus\{j\}}b_i$, $a'':=\bigwedge_{i\in X\setminus\{j\}}b_i\wedge a_j$ and $b:=\bigwedge_{i\in Y}b_i$. 
Then
\begin{eqnarray*}
\mathrm{Pr}(a\wedge b)&=&\mathrm{Pr}(a'\wedge b)-\mathrm{Pr}(a''\wedge b)=\mathrm{Pr}(a')\mathrm{Pr}(b)-\mathrm{Pr}(a'')\mathrm{Pr}(b)\\
&=&(\mathrm{Pr}(a')-\mathrm{Pr}(a''))\mathrm{Pr}(b)=\mathrm{Pr}(a)\mathrm{Pr}(b)
\end{eqnarray*}
as required. The second equality follows by induction. 
\eproof

\noindent
This proposition is in fact a~particular case of Corollary~\ref{cor_corofBITsimpl}.

We arrive at the second main result of our article, the {\em Boolean Independence Theorem}, abbreviated BIT. Recall the evaluation maps $\Phi_A$ of Definition~\ref{def_evalMBool}.

\begin{thm}\label{thm_indepThm}
Let $(A,\mathrm{Pr})$ be a~{\em FPS},
$t$ and $u$ be terms in $T[x_1,\ds,x_n]$ with $n\in\N$ and let $a_1$, $\ds$, 
$a_{2n}$ be events in $A$ such that the events $a_1$, $\ds$, $a_n$ are independent of the events $a_{n+1}$, $\ds$, $a_{2n}$. Then the two events
$$
a:=\Phi_A[t;\,a_1,\,\ds,\,a_n]\,\text{ and }\,b:=\Phi_A[u;\,a_{n+1},\,\ds,\,a_{2n}]
$$ 
are independent.
\end{thm}
\proof
Let $(A,\mathrm{Pr})$, $t$, $u$ and $a_1$, $\ds$, $a_{2n}$ be as stated. 
By Theorem~\ref{thm_DNF} there exist $n$-standard DNF terms $v$ and $w$ such that
$t\stackrel{f}{\approx}v$ and $u\stackrel{f}{\approx}w$. So 
$$
v\in\bigvee_{i=1}^k\mu_i\,\text{ and }\,
w\in\bigvee_{i=1}^l\nu_i
$$
where $k,l\in\N_0$ and the $\mu_i$, resp. $\nu_i$, are 
$[n]$-$\mathrm{monomials}_{\wedge}$ with mutually distinct types. We obtain the term $w'$ by replacing in $w$ every occurrence of the 
variable $x_i$, $i\in[n]$, with the variable $x_{i+n}$. By Proposition~\ref{prop_longDistr} there is a~term $x$ such 
that $(v\wedge w')\stackrel{f}{\approx}x$ and 
$$
x\in\bigvee_{i,\,j=1}^{k,\,l}(\mu_i\wedge\nu_j')
$$
where each $\{n+1,\ds,2n\}$-$\mathrm{monomial}_{\wedge}$ $\nu_j'$ arises from $\nu_j$ by adding $n$ to the indices of the variables $x_1$, $\ds$, $x_n$.

We set $\overline{a}:=(a_1,\ds,a_n)$, $\overline{b}:=(a_{n+1},\ds,a_{2n})$ and $\overline{c}:=(a_1,\ds,a_{2n})$. We see that the events $a$ and $b$ are independent because
\begin{eqnarray*}
\mathrm{Pr}(a\wedge b)&=&   
\mathrm{Pr}\big(\Phi_A[v;\,\overline{a}]\wedge\Phi_A[w;\,\overline{b}]\big)=\mathrm{Pr}\big(\Phi_A[(v\wedge w');\,\overline{c}]\big)=\mathrm{Pr}\big(\Phi_A[x;\,\overline{c}]\big)\\
&=&\sum_{i,\,j=1}^{k,\,l}\mathrm{Pr}\big(\Phi_A[(\mu_i\wedge\nu_j');\,\overline{c}]\big)=\sum_{i,\,j=1}^{k,\,l}\mathrm{Pr}\big(\Phi_A[\mu_i;\,\overline{c}]\big)\mathrm{Pr}\big(\Phi_A[\nu_j';\,\overline{c}]\big)\\
&=&\sum_{i=1}^k\mathrm{Pr}\big(\Phi_A[\mu_i;\,\overline{a}]\big)
\sum_{j=1}^l\mathrm{Pr}\big(\Phi_A[\nu_j;\,\overline{b}]\big)\\
&=&
\mathrm{Pr}\big(\Phi_A[v;\,\overline{a}]\big)\mathrm{Pr}\big(\Phi_A[w;\,\overline{b}]\big)=\mathrm{Pr}(a)\mathrm{Pr}(b)\;.
\end{eqnarray*}
In the first equality we use Proposition~\ref{prop_feqImplseq}. In the 
second equality we use Definition~\ref{def_evalMBool}. In the third equality we use Proposition~\ref{prop_feqImplseq}. In the 
fourth equality we use the fact that $\mu_i\wedge\nu_j'$ are $[2n]$-$\mathrm{monomials}_{\wedge}$ with distinct types, Definition~\ref{def_evalMBool}, Proposition~\ref{prop_diffTypes0} and Corollary~\ref{cor_disjEven}. In the 
fifth equality we use Definitions~\ref{def_evalMBool} and 
\ref{def_Indep}, the assumption on mutual independence and 
Proposition~\ref{prop_indepCompl}. In the 
sixth equality we use the distributive law in the ring $\R$. In the 
seventh equality we use Propositions~\ref{prop_feqImplseq} and \ref{prop_diffTypes0} and Corollary~\ref{cor_disjEven}. Finally, in the eighth equality we use Proposition~\ref{prop_feqImplseq}.
\eproof

We have the following generalization.

\begin{cor}\label{cor_corofBIT}
Let $(A,\mathrm{Pr})$ be a~{\em FPS},
$t_1$, $\ds$, $t_k$ and $u_1$, $\ds$, $u_l$ be terms in $T[x_1,\ds,x_n]$ with $k,l,n\in\N$ and let $a_1$, $\ds$, 
$a_{2n}$ be events in $A$ such that the events $a_1$, $\ds$, $a_n$ are independent of the events $a_{n+1}$, $\ds$, $a_{2n}$. Then the two tuples of events    
$$
\{\Phi_A[t_i;\,a_1,\,\ds,\,a_n]\;|\;i\in[k]\}\,\text{ and }\,
\{\Phi_A[u_i;\,a_{n+1},\,\ds,\,a_{2n}]\;|\;i\in[l]\}
$$
are independent one of another.
\end{cor}
\proof
We set $b_i:=\Phi_A[t_i;a_1,\ds,a_n]$, $i\in[k]$, and 
$c_i:=\Phi_A[u_i;a_{n+1},\ds,a_{2n}]$, $i\in[l]$. For any subsets $X\sus[k]$ 
and $Y\sus[l]$, the two events
$$
{\textstyle
\bigwedge_{i\in X}b_i=\Phi_A\big[\bigwedge_{i\in X}t_i;\,a_1,\,\ds,\,a_n\big]
}
$$
and 
$$
{\textstyle
\bigwedge_{i\in Y}c_i=\Phi_A\big[\bigwedge_{i\in Y}u_i;\,a_{n+1},\,\ds,\,a_{2n}\big]
}
$$
are independent by the previous theorem. Here we represent the argument 
$\bigwedge_{i\in X}t_i$ and $\bigwedge_{i\in Y}u_i$ by any term in it. This is possible by Propositions~\ref{prop_feqImplseq} and \ref{prop_longConj}
\eproof

\noindent
We state this corollary equivalently and less formally.

\begin{cor}\label{cor_corofBITsimpl}
If $n\in\N$ and some events $a_1$, $\ds$, $a_n$ in a~{\em FPS} $(A,\mathrm{Pr})$ are independent of 
some other events $a_{n+1}$, $\ds$, $a_{2n}$ and 
$$
A_1:=\langle\{a_1,\,\ds,\,a_n\}\rangle\,\text{ and }\,A_2:=\langle\{a_{n+1},\,\ds,\,a_{2n}\}\rangle
$$ 
are the corresponding (finite) Boolean subalgebras of $A$ generated by the two 
$n$-tuples, then any tuple of events $X\sus A_1$ is independent of any tuple $Y\sus A_2$.
\end{cor}

\section[The BIT in the Probabilistic Method]{The BIT in the Probabilistic Method}\label{sec_PM}

The crux of independence of events is the following simple and well known 
enumerative proposition. We omit the proof.

\begin{prop}\label{prop_keyIndep}
If $X$ and $Y$ are finite sets then $|X\times Y|=|X|\cdot|Y|$.    
\end{prop}

We give just one application of 
Corollary~\ref{cor_corofBITsimpl} in the Probabilistic Method, but more can be mined from \cite{alon_spen}. Following 
\cite[Theorem~5.2.1]{alon_spen} and the notation there, we consider a~hypergraph $H=(V,E)$, which 
is a~pair of finite sets $V\ne\emptyset$ and $E\sus\mathcal{P}(V)$, such that every 
edge $f\in E$ has $|f|\ge k$ elements and intersects at most $d$ other edges. One then proves by means of LLL (the Lov\'asz Local Lemma) 
that if $\mathrm{e}(d+1)\le 2^{k-1}$ (where $\mathrm{e}=2.71828\ds$ is the 
Euler number) then there exists a~map (a~vertex coloring)
$$
\chi\cc V\to\{\mathrm{r},\,\mathrm{b}\}=:C
$$
that is non-constant on every edge $f\in E$.

We refer to \cite{alon_spen} for the precise statement of LLL; see \cite[Chapter~6]{mitz_upfa} for applications of LLL in Computer Science.   The independence of an event of a~tuple of events, which 
we reviewed earlier (and which is undefined in \cite{alon_spen}), plays key role in it. We quote from the proof in \cite{alon_spen}:  ``Color each vertex $v$ of $H$, randomly and independently, either blue 
or red (with equal probability). For each edge $f\in E$, let $A_f$ be the event 
that $f$ is monochromatic.'' (i.e., the coloring is constant on $f$) $\ds$ 
``Moreover, each event $A_f$ is clearly mutually independent of all the other 
events $A_{f'}$ for all edges $f'$ that do not intersect $f$.''

We make these arguments more formal (i.e., more precise) in order to see how they fit in
Corollary~\ref{cor_corofBITsimpl}; we in fact formally prove, probably for the 
first time ever in the history of applications of LLL, that the event $A_f$ 
is indeed independent of  all events $A_{f'}$ with $f'\cap f=\emptyset$ and 
therefore the hypothesis of LLL is satisfied.  We set $|V|=n\in\N$ and use 
a~FPS $(\mathcal{P}(M),\mathrm{Pr})$ of the kind reviewed earlier 
($C=\{\mathrm{r},\mathrm{b}\}$):
$$
M:=\{\chi\;|\;\chi\cc V\to C\}\,\text{ and always }\,\mathrm{Pr}(\{\chi\}):=2^{-n}\;.
$$
For every $v\in V$ and 
$c\in C$ we denote by $A_{v,c}$ the event that $\chi(v)=c$,
$$
A_{v,\,c}=\{\chi\in M\;|\;\chi(v)=c\}\;.
$$
We claim that if $X,Y\sus V$ are disjoint sets, then the two tuples of events
$$
B_X:=\{A_{v,\,c}\;|\;v\in X,\,c\in C\}\,\text{ and }\,B_Y=\{A_{v,\,c}\;|\;v\in Y,\,c\in C\}
$$
are independent one of another. We postpone the proof of it for 
a~moment and deduce at once, using Corollary~\ref{cor_corofBITsimpl}, that for any $f\in E$ the 
event
$$
{\textstyle
A_f=\bigwedge_{v\in f}A_{v,\,\mathrm{r}}\vee\bigwedge_{v\in f}A_{v,\,\mathrm{b}}
}
$$
is indeed independent of the tuple of events
$$
{\textstyle
\{A_{f'}\;|\;f'\in E,\,f\cap f'=\emptyset\}\;.
}
$$

We show that the tuple $B_X$ is independent of $B_Y$. For any subsets $U\sus X\times 
C$ and $W\sus Y\times C$ we set
$$
M_U:=\{\chi\in M\;|\;(v,\,c)\in U\Rightarrow\chi(v)=c\}
$$
and define similarly $M_W$ and $M_{U\cup W}$. We also set
$$
N_U:=\{\chi\cc X\to C\;|\;(v,\,c)\in U\Rightarrow\chi(v)=c\}
$$
and define similarly $N_W$ and $N_{U\cup W}$ (as sets of maps with respective definition domains $Y$ and $X\cup Y$).
Then
\begin{eqnarray*}
&&{\textstyle
\mathrm{Pr}\big(\bigwedge_{(v,\,c)\in U}A_{v,\,c}\wedge\bigwedge_{(v,\,c)\in W}A_{v,\,c}\big)
}=\frac{|M_{U\cup W}|}{2^n}
=\frac{|N_{U\cup W}|\cdot2^{|V\setminus(X\cup Y)|}}{2^n}\\
&&=\frac{|N_{U\cup W}|}{2^{|X\cup Y|}}=
\frac{|N_U|}{2^{|X|}}\cdot
\frac{|N_W|}{2^{|Y|}}=\frac{|N_U|\cdot 2^{|V\setminus X|}}{2^{n}}\cdot
\frac{|N_W|\cdot 2^{|V\setminus Y|}}{2^n}\\
&&=\frac{|M_U|}{2^n}\cdot\frac{|M_W|}{2^n}=
{\textstyle\mathrm{Pr}\big(\bigwedge_{(v,\,c)\in U}A_{v,\,c}\big)\mathrm{Pr}\big(\bigwedge_{(v,\,c)\in W}A_{v,\,c}\big)}
\end{eqnarray*}
and $B_X$ is independent of $B_Y$. 

The key fourth equality in the computation holds due to $X\cap Y=\emptyset$. Then 
$2^{|X\cup Y|}=2^{|X|}\cdot2^{|Y|}$ and $|N_{U\cup W}|=|N_U|\cdot|N_W|$. 
This last equality immediately follows from  Proposition~\ref{prop_keyIndep}.

\section{Concluding remarks}\label{sec_conclCom}

We hope to extend the theory of formal polynomials in 
Sections~\ref{sec_defini}--\ref{sec_polynomials} to rings of formal 
and quasi-formal power series.

In the independence part in 
Sections~\ref{sec_BooTermPoly}--\ref{sec_PM} we answered to our satisfaction the question what 
is the precise role and form of independence of events in the 
Probabilistic Method, especially in LLL and its applications. It was important for 
us since this author belongs to those mathematicians who care more about the 
proof of a~theorem or the reason why a~computation works than about the 
claim of the theorem or the result of the computation. 

In the literature we did not encounter any signs of dissatisfaction with LLL and its 
applications, in the sense that they are (in the present form) non-rigorous; but 
they clearly are. It is an interesting psychological question why it is so 
obvious to us, the users of LLL, that the (above) event $A_f$ corresponding to an edge $f$ 
is independent of all events $A_{f'}$ with $f'$ disjoint to $f$. It seems that we employ some kind of mental alliteration or rhyme: disjoint --- independent. Rhymes 
sound good but of course do not prove anything. Try a~psycho-linguistic 
experiment: instead of the relational adjective ``independent of'', name the 
independence of events, $\mathrm{Pr}(a\wedge b)=\mathrm{Pr}(a)\mathrm{Pr}(b)$ 
and so on, by a~nonsensical term, like ``knowledgeable of'', or even by the opposite 
``dependent on''. Is it then so clear that if $f$ and $f'$ are disjoint edges then the 
event $A_f$ is knowledgeable of $A_{f'}$,
or that $A_f$ is dependent on $A_{f'}$? In the near future we hope to rigorize or make 
more precise arguments and proofs for some other results in the Probabilistic Method.

\bigskip\noindent
{\em Department of Applied Mathematics\\ 
Faculty of Mathematics and Physics\\
Charles University\\ 
Malostransk\'e n\'am\v est\'\i\ 25\\
118 00 Praha\\
Czechia}\\
{\tt klazar@kam.mff.cuni.cz}
\end{document}